\documentclass[11pt]{article}
\usepackage[T1]{fontenc}

\usepackage[utf8]{inputenc}

\usepackage{amsthm}

\usepackage{mathrsfs}
\usepackage{amsmath}
\usepackage{amssymb}
\usepackage{graphicx}
\usepackage{color}
\usepackage{caption,subcaption}
\usepackage{babel}
\usepackage{comment}
\usepackage{algorithm}
\usepackage{algpseudocode}
\usepackage{setspace}
\usepackage{grffile}
\usepackage{enumerate}

\usepackage{hyperref}
\hypersetup{colorlinks=true,
            linkcolor=black,
            anchorcolor=black,
            citecolor=black}

\numberwithin{equation}{section}
\numberwithin{figure}{section}
\numberwithin{table}{section}

\newtheorem{THEOREM}{Theorem}[section]

\newtheorem{proposition}[THEOREM]{Proposition}

\newtheorem{remark}[THEOREM]{Remark}

\newcommand{\R}{\mathbb{R}}

\newcommand{\lam}{\lambda}

\newcommand{\bx}{\mathbf{x}}

\newcommand{\bu}{\mathbf{u}}

\newcommand{\sig}{\sigma}
\newcommand{\gam}{\gamma}
\newcommand{\bv}{\mathbf{v}}
\renewcommand{\th}{\theta}

\newcommand{\lan}{\langle}
\newcommand{\ran}{\rangle}
\newcommand{\eps}{\varepsilon}

\newcommand{\bn}{\mathbf{n}}

\newcommand{\Dx}{\Delta x}
\newcommand{\Dt}{\Delta t}
\newcommand{\Dv}{\Delta v}

\newcommand{\rd}{\mathrm{d}}

\newcommand{\ion}{\rho_{ion}}
\newcommand{\hphi}{\widehat{\phi}}

\newcommand{\hh}{\widehat{h}}

\newcommand{\cH}{\mathcal{H}}
\newcommand{\cJ}{\mathcal{J}}
\newcommand{\cL}{\mathcal{L}}

\newcommand{\feq}{\overline{f}}
\newcommand{\balpha}{\mathbf{\alpha}}

\newcommand{\dE}{\delta E}

\newcommand{\df}{\delta f}
\newcommand{\dfd}{\delta f^\Delta}
\newcommand{\drho}{\delta \rho}
\newcommand{\dH}{\delta H}

\newcommand{\baE}{\overline{E}}

\newcommand{\norm}[1]{\left\|#1\right\|}
\newcommand{\inner}[1]{\left\langle #1\right\rangle}

\newcommand{\email}[1]{\href{mailto:#1}{\nolinkurl{#1}}}

\title{Dynamical feedback control with operator learning for the Vlasov-Poisson system}

\author{Jingcheng Lu\thanks{School of Mathematics, University of Minnesota–Twin Cities, MN;
Institute of Mathematical Sciences, ShanghaiTech University, China
  (\email{lu000688@umn.edu}).}
\and Li Wang \thanks{School of Mathematics, University of Minnesota–Twin Cities, MN
  (\email{liwang@umn.edu}).}
\and Jeff Calder\thanks{School of Mathematics, University of Minnesota–Twin Cities, MN
  (\email{jwcalder@umn.edu}).}}
\date{}

\begin{document}

\maketitle

\begin{abstract}

To meet the demands of instantaneous control of instabilities over long time horizons in plasma fusion, we design a dynamic feedback control strategy for the Vlasov–Poisson system by constructing an operator that maps state perturbations to an external control field. In the first part of the paper, we propose learning such an operator using a neural network. Inspired by optimal control theory for linearized dynamics, we introduce a low-rank neural operator architecture and train it via adjoint state method. The resulting controller is effective at suppressing instabilities well beyond the training time horizon.  To
generalize control across varying initial data, we further introduce a novel cancellation-based control strategy that removes the destabilizing component of the electric field. This approach naturally defines an operator without requiring any training, ensures perturbation decay over infinite time, and demonstrates strong robustness under noisy feedback. Numerical experiments confirm the effectiveness of the method in both one- and multidimensional settings.

\end{abstract}

\tableofcontents

\section{Introduction}
Long-term suppression of plasma instabilities remains a major challenge in achieving controlled nuclear fusion. Without adequate confinement, the plasma can gradually deviate from its desired state, potentially leading to a breakdown of the reaction or a significant reduction in fusion efficiency. For example, the two-stream instability can cause unwanted beam scattering and loss of focus \cite{sydorenko2016effect}. Likewise, the bump-on-tail instability, often driven by runaway electrons or radiofrequency heating, can reduce plasma heating efficiency, degrade confinement, and cause significant energy deposition on the reactor walls. Therefore, it is crucial to control such instabilities through the application of carefully designed external forces.

To place the problem on a more concrete footing, we consider the Vlasov–Poisson (VP) system, which is fundamental for describing collisionless plasma dynamics:
\begin{equation}\label{eq:vlasov poisson}
\left\{
\begin{array}{l}
\partial_t f(\bx,\bv,t)+\bv\nabla_\bx f(\bx,\bv,t)+(E(\bx,t)+H(\bx,t))\cdot\nabla_\bv f(\bx,\bv,t) = 0 \,,\\
\\
E(\bx,t) = -\nabla_\bx \Phi(x,t), \quad \nabla_\bx E(\bx,t) = -\Delta\Phi(\bx,t) = \rho(\bx,t)-\ion\,,\\
\\
\rho(\bx,t) = \int f(\bx,\bv,t) d\bv\,.
\end{array}
\right.
\end{equation}
Here $f(\bx,\bv,t)$ is particle density at location $\bx$, time $t$, with velocity $\bv$, $E(\bx,t)$ is the self-generated electric field, $\rho(\bx,t)$ is the charge density, and $\ion$ is the constant density of background ion. $H(\bx,t)$ is the external electric field that will be used for the control purposes. 

For ease of discussion, we assume that the desired confined plasma state is described by the distribution $\feq(\bx,\bv)$, which is an equilibrium of the system \eqref{eq:vlasov poisson} when $H = 0$. There has been an extensive study on the stability or instability of the equilibrium state  $\feq(\bx,\bv)$  upon perturbation. In particular, when $\feq(\bx, \bv)$ depends only on $\bv$ and follows a Maxwellian distribution, i.e., $\feq(\bx, \bv)= \feq(\bv)\propto e^{-{|\bv-u|^2/2T}}$ for some given bulk velocity $u$ and temperature $T$, this equilibrium is stable, in the sense that small perturbations in plasma beams lead to the damping of electrostatic waves, a phenomenon known as the Landau damping \cite{Landau, mouhot2011landau}. In contrast, when $\feq(\bx, \bv) = \feq(\bv)$ is a mixture of two Gaussians, the equilibrium becomes unstable; that is, an initial disturbance in the beams is amplified, leading to undesirable growth of the perturbation and the transfer of energy from the plasma beams to electrostatic waves \cite{chen2016}. 
Additionally, \eqref{eq:vlasov poisson} can also admit spatially inhomogeneous equilibrium, which may be unstable \cite{guo1995instability} under certain conditions.

 To suppress the growth of perturbation over a long time horizon $t\in[0,T]$, we aim to design an external electric field $H(x,t)$ via the PDE-constrained optimization framework:
 \begin{equation}\label{eq:VP constrained optim}
 \min_{H\in L^2_{\bx,t}}  \ \cJ(T;H)
 \quad \text{s.t.} \quad f ~\text{solves}~~ \eqref{eq:vlasov poisson}\,.
 \end{equation}
The loss functional $\cJ(T;H)$ defines the control objective, which may include terms such as a running loss $\displaystyle \frac{1}{2}\int^T_0 ||f(t;H)-\feq||^2_{\bx,\bv} \rd t$, a terminal loss $\displaystyle \frac{1}{2}||f(T;H)-\feq||^2_{\bx,\bv}$, and potentially others. A recent study has analyzed how the choice of objective influences the optimization landscape \cite{guerra2025metric}.
The parameterization of $H$ is also crucial for effective control. Directly assuming $H$ to be an arbitrary function of $x$ (or of $x$ and $t$) typically leads to optimization landscapes with many local minima. As observed in \cite{einkemmer2024suppressing}, a parameterization of the form
\begin{align} \label{H0}
    \displaystyle H(\bx,t;\balpha) = \sum^r_{k=1} \alpha_k\phi_k(\bx,t)
\end{align}
with $\alpha_k$ as the parameters to be optimized, not only reduces the dimensionality of the control problem but also tends to convexify the optimization landscape. A moment-based formulation was also proposed in \cite{lu2025controlling}, offering a potential reduction of the dimension of \eqref{eq:VP constrained optim}. 

Despite the above efforts in selecting the objective function $\cJ$ for the parameterization of $H$, a key drawback remains: the optimized field $H$ is only effective over the time interval $[0,T]$. To extend its usefulness beyond $T$, one would need to solve \eqref{eq:VP constrained optim} again, which is computationally expensive and incompatible with the practical requirement of near-instantaneous control. For this reason, dynamical feedback control appears to be a more practical alternative, provided such a control can be constructed. Several studies have investigated this direction. One approach is to determine a control $H(t,\bx)$ over short intervals $[t,t+\Delta t]$, where the objective function is restricted to this local time window, instead of seeking a control valid over the entire time span $[0,T]$. By discretizing the problem and solving the corresponding optimality conditions, an exact relation between $H$ and the solution $f$ can be derived, which then serves as the feedback law \cite{albi2025instantaneous}.
Another line of work employs reinforcement learning, which assumes the underlying model is unknown and instead learns the feedback law in a data-driven manner \cite{dubbioso2023deep, seo2024avoiding}. 
Similar model-free strategies have also been applied to tokamak magnetic control and plasma-shaping tasks \cite{tracey2024towards,degrave2022magnetic,char2023offline}.

In this paper, we propose a novel approach to deriving the feedback control law by designing and learning a neural operator. Specifically, we aim to learn a mapping from deviation of the particle distribution from equilibrium $\delta f(t,\bx, \bv):= f(t,\bx,\bv)- \feq(t,\bx,\bv)$ to $H(t,\bx)$ such that the control objective $\cJ$ is effectively minimized. This constitutes a mapping between function spaces: from the space in which $\delta f$ resides to the space where $H$ is defined, i.e., $H[\df](t,\bx)$. A natural way to approximate such a mapping is through a \emph{neural operator}, which refers to a neural network that is designed to learn a mapping between function spaces, i.e., an operator, instead of a function between Euclidean spaces. Neural operators have been extensively studied in recent years across a wide range of applications. Several architectures have been proposed, including the Fourier Neural Operator (FNO) \cite{li2020fourier,kovachki2021universal,li2023fourier}, Low-Rank Neural Operator (LNO) \cite{kovachki2023neural}, Graph Neural Operator (GNO) \cite{anandkumar2020neural,li2020multipole}, and DeepONet \cite{lu2021learning,wang2021learning,lanthaler2022error}. Comprehensive overviews can be found in two recent review articles \cite{lu2022comprehensive, kovachki2023neural}.

To construct such an operator for our setting, we design a specific architecture guided by linear  optimal control  theory. In particular, we first derive a feedback control law for the linearized system, and then use it as a blueprint to design a neural operator for the general nonlinear problem. The neural operator is then trained via an adjoint-based method. The resulting learned neural operator remains effective over much longer time horizons than directly learning the control field $H$ without feedback information, and it is more robust to noise. Finally, to extend the control performance to arbitrarily long times and arbitrary initial perturbations, we introduce a quasi-optimal law based on a cancellation principle.

The rest of the paper is organized as follows. In the next section, we construct a neural operator to obtain the feedback law, with its structure inspired by linear optimal control theory and its implementation carried out through an adjoint-based method. The effectiveness of the resulting neural operator is demonstrated through extensive examples in Section~\ref{sec:numerical exp}. In Section~\ref{sec4}, we propose a novel universal feedback law based on a cancellation argument. This law can be viewed as a special case of a neural operator and is training-free. Numerical experiments further confirm the effectiveness of this approach.
 

\section{Dynamical feedback control with low-rank neural operator}
In this section, we investigate the structure of the feedback control $\df \mapsto H$. For clarity of presentation, we focus on the one-dimensional case, and assume a spatially homogeneous equilibrium, i.e., $\feq(x,v) \equiv \feq(v)$, although the derivations extend naturally to higher dimensions.

\subsection{Operator-based external field}
To gain primary insight about the appropriate structure of the mapping from deviation $\df$ to the control field $H(\df(t))(t,x)$:
\[
H(\cdot )(\cdot): ~ \df \mapsto H\,,
\]
we first consider the following linearized Vlasov–Poisson system around the target equilibrium $\feq(v)$. More specifically, write $f(t,x,v) = \feq(v) + \df(t,x,v)$, and assume that $\|\df\| \ll \| \feq \|$, we have the following linearized equation for the perturbation $\df$:
\begin{equation}\label{eq:vlasov poisson linear}
\left\{
\begin{array}{l}
\partial_t \df+v\partial_x \df+(\dE+H)\partial_v \feq = 0 \,,\\
\\
\dE(x,t) = -\partial_x \Phi(x,t), \quad -\partial^2_x\Phi(x,t) = \drho(x,t) \,,\\
\\
\drho(x,t) = \int^\infty_{-\infty} \df(x,v,t) dv\,.
\end{array}
\right.
\end{equation}
Assuming the control objective takes the form:
\begin{equation}\label{eq:quadratic loss}
\cJ(\df,H) = \int^\infty_0 \frac{1}{2}||\df(\cdot,\cdot, t)||^2_{x,v}+
\frac{\beta}{2}||H(\cdot,t)||^2_x \rd t,
\end{equation}
where $\beta>0$ is the weight that penalizes the magnitude of the control $||H||$. We can define the corresponding Hamiltonian: 
\begin{equation}\label{eq:hamiltonian}
\begin{split}
\cH(\df,H,\lam) = 
-\iint\frac{1}{2}\df(x,v,t)^2\rd xdv& 
-\int \frac{\beta}{2}H(x,t)^2 \rd x \\
&-\iint \lambda(x,v,t)\big(v\partial_x\df+(\dE+H)\partial_v\feq\big)\rd x dv.
\end{split}
\end{equation}
Let $G(x_1,x_2)$ be the Green's function of the Poisson equation with suitable boundary conditions, the perturbation of the self-generated field $\dE$ can be written as
\[
\dE(x,t) = \int \partial_{x_1}G(x,y)\delta\rho(y,t)\rd y = \iint \partial_{x_1} G(x,y)\df(y,v',t)\rd y \rd v'\,.
\]
By interchanging of variables $(x,v)$ and $(y,v')$, we obtain
\[
\iint \lambda(x,v,t)\dE(x,t)\partial_v\feq(v)\rd x \rd v = \iint \big[\iint \partial_{x_1} G(y,x)\lambda(y,v',t)\partial_v\feq(v')\rd y dv'\big]\df(x,v,t)\rd x \rd v.
\]
Then the Hamiltonian \eqref{eq:hamiltonian} can be rewritten as 
\[
\begin{split}
\cH(\df,H,\lam) &= 
-\iint\frac{1}{2}\df(x,v,t)^2\rd x \rd v-\int \frac{\beta}{2}H(x,t)^2 \rd x \\
&+\iint \big[v\partial_x\lambda(x,v,t)-\lan\partial_{x_1}G(\cdot,x),\int \lambda(\cdot,v',t)\partial_v\feq(v')dv' \ran\big]\df(x,v,t)\rd x  \rd v\\
&-\iint \lambda(x,v,t)H(x,t)\partial_v\feq(v)\rd x \rd v.
\end{split}
\]
By the Pontryagin Maximum Principle, the optimal control, $H^*$, and the corresponding adjoint variable, $\lambda^*$, satisfy the conditions $\displaystyle \partial_t \lambda = -\frac{\partial \cH}{\partial \df}$ and $\displaystyle \frac{\partial \cH}{\partial H} = 0$, that is,
\begin{equation}\label{eq:linearized eqn optim cond}
\begin{split}
\left\{\begin{array}{l}
\partial_t \lambda^*+v\partial_x\lambda^*-\lan\partial_{x_1}G(\cdot,x),\int \lambda^*(\cdot,v',t)\partial_v\feq(v')\rd v '\ran=
\df\,,  \\
\\
-\beta H^*-\int\lam^*\partial_v\feq dv= 0\,.
\end{array}
\right.
\end{split}
\end{equation}
The first equation shows the linear dependence of $\lambda^*$ on $\df$, and the second equation indicates the linear dependence of $H^*$ on $\lambda^*$. Consequently, the optimal control $H^*$ for the linearized system \eqref{eq:vlasov poisson linear} is a linear functional of the perturbation $\df$, although an explicit representation of the mapping $\df \mapsto H^*$ is not readily available. 

Nevertheless, this linear dependence motivates us to construct 
the feedback control mechanism through a single-layer low-rank neural operator \cite{kovachki2023neural}:
\begin{equation}\label{eq:H LNO}
H[\df(t)](x;\th) = \sum^{r}_{k=1}\phi_k(x;\th_\phi) \iint \psi_k(y,v;\th_\psi)\df(y,v,t) \rd y \rd v, \quad \th = \{\th_\phi,\th_\psi\}.
\end{equation}
Here, $\phi^{NN}(x; \theta_\phi) := \big(\phi_1(x; \theta_\phi), \phi_2(x; \theta_\phi), \ldots, \phi_r(x; \theta_\phi)\big)$ is a learnable mapping $\R$ to $\R^r$ consisting of basis functions parameterized by $\theta_\phi$, while the kernel $\psi^{NN}(x,v; \theta_\psi) = \big(\psi_1(x,v; \theta_\psi), \psi_2(x,v; \theta_\psi), \ldots, \psi_r(x,v; \theta_\psi)\big)$ is a learnable mapping from $\R^2$ to $\R^r$, parameterized by $\theta_\psi$. In the following, we will parameterize $\phi^{NN}$ and $\psi^{NN}$ by neural networks. When considering a bounded domain $x \in \Omega$ with specific boundary conditions, it is preferable for the basis $\phi^{NN}(x; \theta_\phi)$ to explicitly encode these boundary conditions. We emphasize that, unlike \eqref{H0}, the form \eqref{eq:H LNO} introduces time dependence in 
$H$ through $\df$. This provides a more effective way to incorporate time dependence into the control than simply allowing $\alpha$ to vary with time in \eqref{H0}. The latter approach not only increases the dimensionality of the optimization problem, since discretizing time into $N_t$ steps enlarges the parameter space by a factor of $N_t$, but may also fail to capture meaningful time dependence, as direct optimization can yield coefficients that vary little with $t$.

In practical computations, we have perturbations only on the computational grid: 
\[
\dfd = \{\df(x_i,v_j)\}_{1\leq i\leq N_x, 1\leq j\leq N_v}\in \R^{N_x N_v}.
\]
Then the nodal values of the external field $\displaystyle H^\Delta = \{H[\dfd](x_i)\}_{1\leq i \leq N_x}$ can be represented by:
\[
H^\Delta_i = \sum^r_{k=1}\phi_k(x_i)\sum_{\substack{1\leq l\leq N_x
\\ 1\leq m\leq N_v}}\psi_k(x_l,v_m)\df(x_l,v_m)\Dx\Dv, \quad 1\leq i \leq N_x\,.
\]
Equivalently, this can be written in the matrix form as
\begin{equation}\label{eq:H low rank matrix}
H^\Delta = \left([\phi^\Delta_{r\times N_x}]^\top \psi^\Delta_{r\times N_xN_v}\Dx\Dv\right)\dfd,
\end{equation}
where $\phi^\Delta_{r\times N_x} = [\phi_k(x_i)]_{k,i}$ and $\psi^\Delta_{r\times N_x N_v} = [\psi_k(x_i,v_j)]_{k,i,j}$ are two projection matrices. 

Note that \eqref{eq:H low rank matrix} is reminiscent of the well-known Linear Quadratic Regulator (LQR) \cite{kalman1960new,kalman1961new} for linear ODE systems, where the optimal control problem is given by
\[
 \begin{split}
 &\min_{a(t)}  \int^\infty_0 g(t)^\top Q g(t)+a(t)^\top R a(t) \rd t\,, \\
s.t. & \qquad \qquad \frac{dg}{\rd t} = Ag(t)+Ba(t)\,,
 \end{split}
\]
where $Q$ and $R$ are positive definite matrices. 
If we assume that the optimal control takes the form  $a(t) = -Kg(t)$, then the feedback matrix $K$ can be solved from the algebraic Riccati equation.
In the context of the semi-discrete ODE arising from the linearized equation \eqref{eq:vlasov poisson linear}, the controller \eqref{eq:H low rank matrix} corresponds to a linear projection from $\dfd\in\R^{N_xN_v}$ to $H^\Delta\in \R^{N_x}$, constructed using low-rank matrices.

We emphasize that the choice of the low-rank linear neural operator architecture is inspired by both mathematical and physical considerations. We summarize the key components:

\medskip
\noindent \textbf{1. Low-rank/finite-basis structure.} A key aspect of our operator design is to represent the feedback map in the form
\[
H[\delta f(t)](x) =\sum_k \alpha_k(\delta f(t))\phi_k(x).
\]
Here, the basis functions $\phi_k(x)$ represent the spatial modes available to the control field. This is useful both mathematically and practically: the choice of $\phi_k$ can reflect boundary conditions, and in applications, the imposed field may only be generated through a finite number of modes or actuator channels. 

This form is also related to the structure of DeepONet-type representations,
\[
H[\delta f(t)](x)=\sum_k t_k(x)b_k(\delta f(t)),
\]
where the trunk functions $t_k(x)$ correspond to the output basis functions $\phi_k(x)$. In our construction, the branch coefficients $b_k(\delta f)$ are not represented by a generic nonlinear network, but are obtained through integral kernels applied to the $(x,v)$-dependent perturbation. This makes the information extraction step directly tied to the kinetic state and preserves an interpretable feedback structure.

\medskip
\noindent \textbf{2. Extract information via the integral kernel.} The controller \eqref{eq:H LNO} uses integral kernels to extract the coefficients $\alpha_k(\delta f)$. The controller maps a kinetic perturbation $\delta f(x,v)$, defined over position and velocity variables, to a spatial control field $H(x)$. This dimension-reducing structure is not the native setting of many standard operator architectures, such as U-Net or FNO, which are typically designed to map between functions defined on equal dimensions. 

By contrast, the integral-kernel representation provides a direct and interpretable way to perform this reduction: it extracts observables from the $(x,v)$-dependent perturbation and uses them to activate spatial control modes.  It also has potential advantages for nonuniform grids or more complex computational geometries, since the integrals can, in principle, be approximated using appropriate quadrature rules. Then the output basis can be chosen to match the spatial domain. 

\medskip
\noindent \textbf{3. Linear structure.} The linearity of the feedback operator is another essential design
choice, which is first motivated by the linearized optimal control problem. Although the full Vlasov--Poisson system is nonlinear with respect to $\delta f$, the numerical experiments in Section~\ref{sec:numerical exp} show that a
properly optimized linear controller can suppress
perturbations effectively over extended time horizons. Moreover, the
cancellation-based controller \eqref{eq:H cancellation based}, derived in Section~\ref{sec4} from the nonlinear perturbation equation, belongs to the same linear operator class \eqref{eq:H LNO} and satisfies a provable dissipation estimate \eqref{decay-df}. Together, these observations support the use of a linear feedback architecture beyond the linearized derivation.



\subsection{Optimization with adjoint state method}\label{sec:adjoint method}

Using the form \eqref{eq:H LNO} for our neural operator that maps $\df$ to $H$, we need to define a loss function to determine the optimal basis and kernels so that the control remains effective over an extended period and for a wide range of perturbations. Ideally, the control objective should consider a sufficiently diverse family of initial perturbations, such that the perturbations evolving over time remain within the subspace generated by this initial family. At the same time, the control horizon should be long enough to capture and suppress any instabilities that may arise. In practice, this is achieved by adding noise to a given initial condition, as will be explained in detail in Section~\ref{sec:computational setups}. For simplicity of presentation, here we consider the control problem for a single initial condition and define the loss function over a sufficiently long time interval $[0, T]$. 

Specifically, we optimize the external field operator $H_\th[\df]$, through the following problem,
\begin{subequations}\label{eq:fb optim}
 \begin{equation}\label{eq:fb optim loss}
\min_{\th} \quad \cJ(\th) = \frac{1}{2}\int^T_0||f(\cdot,\cdot,t;H_\th)-\feq(\cdot,\cdot)||_{x,v}^2\rd t
\end{equation}
\begin{equation}\label{eq:fb optim constraint}
s.t. \quad  \left\{
\begin{array}{l}
\partial_t f +v\partial_x f +(E[f]+H_\th[\df])\partial_v f = 0\,, \\
\\
E = -\partial_x \Phi, \quad -\partial^2_x\Phi = \rho-\ion \,,\\
\\
\rho = \int f \rd v \,.
\end{array}
\right.
 \end{equation}
\end{subequations}
Note that the running loss \eqref{eq:fb optim loss} corresponds to a finite-time truncation of the quadratic loss \eqref{eq:quadratic loss} with $\beta=0$. In our experiments, we observed that including the control penalty $\displaystyle \frac{\beta}{2} ||H||^2_2$ does not lead to any improvement in computational performance.

To compute the gradient $\nabla_\th \cJ(\th)$, we introduce the Lagrangian multiplier,
\[
\begin{split}
\cL(\th;\lambda) = \frac{1}{2}\int^T_0||f(t;H_\th[\df])-\feq||_{x,v}^2\rd t+\int^T_0 \iint \lam(x,v,t)\big(\partial_t f+v\partial_x f+(E[f]+H_\th[\df])\partial_v f\big)\rd x \rd v \rd t\,.
\end{split}
\]
Note that when $f$ satisfies the constraint PDE \eqref{eq:fb optim constraint}, we have $\cL(\th;\lambda) = \cJ(\th)$, and consequently, $\nabla_\th \cL = \nabla_\th \cJ$, for arbitrary $\lambda$.  Integration by parts yields
\[
\begin{split}
\nabla_\th \cJ = \nabla_\th \cL &= \int^T_0\iint (f(x,v,t)-\feq(x,v))\nabla_\th f(x,v,t) \rd x \rd v \rd t+\iint \lambda(x,v,t)\nabla_\th f(x,v,t)|^T_{t=0}\rd x \rd v \\
&-\int^T_0\iint\big(\partial_t\lambda+v\partial_x\lambda+(E[f]+H_\th[\df])\partial_v\lambda\big)\nabla_\th f\rd x \rd v \rd t\\
& +\int^T_0 \iint \big\{\nabla_\th E(x,t) +H[\nabla_\th \df(t)](x;\th)+(\nabla_\th H)[\df(t)](x;\th)\big\}\lambda(x,v,t)\partial_v f(x,v,t) \rd x \rd v \rd t.
\end{split}
\]
By assuming that  $\lambda$ solves the adjoint  equation
\begin{equation}\label{eq:adjoint pde}
\partial_t \lambda+v\partial_x \lambda+(E[f]+H_\theta[\df])\partial_v \lambda = \df, \quad \lambda(x,v,T) = 0,    
\end{equation}
the gradient of the loss function reduces to
\begin{equation} \label{grad}
\begin{split}
\nabla_\th \cJ = \int^T_0 \iint\big\{\nabla_\th E(x,t) &+H_\theta[\nabla_\th \df(t)](x;\th)\\
&+(\nabla_\th H_\theta)[\df](x;\th)\big\}\lambda(x,v,t)\partial_v f(x,v,t) \rd x \rd v \rd t.
\end{split}
\end{equation}

Note that the adjoint equation \eqref{eq:adjoint pde} and the corresponding gradient \eqref{grad} may appear slightly different from those obtained using the classical adjoint state method. In the classical approach as is done in \cite{einkemmer2024suppressing}, 
one regards
$E$ as a functional of the charge density $\rho$. Therefore,
\[
\nabla_\theta E
=
\frac{\delta E}{\delta \rho}[\nabla_\theta \rho],
\qquad
\nabla_\theta \rho=\int \nabla_\theta f\,dv,
\]
where $\frac{\delta E}{\delta \rho}[\cdot]:L_x^2\to L_x^2$ is a linear operator. After integration by parts in the
adjoint calculation, the adjoint operator
$\left(\frac{\delta E}{\delta \rho}\right)^*$
acts on the corresponding source term and enters the adjoint equation \eqref{eq:adjoint pde}. In this way, $\nabla_\theta E$ itself does not appear explicitly in the final gradient formula \eqref{grad}.

Likewise, the term $H_\theta[\nabla_\th \df(t)](x;\th)$ can be eliminated from \eqref{grad} by explicitly writing: 
\begin{align*}
H_\theta[\df](x;\th) &= \sum_r \phi_r(x)\iint \psi_r(y,v)\df(y,v)\rd y \rd v
\\ & = \sum_r \phi_r(x)\iint \psi_r(y,v)(f(y, v) - \feq(v))\rd y \rd v
\end{align*}
and then 
\begin{align*}
H_\theta[ \nabla_\theta \df](x;\th) = \sum_r \phi_r(x)\iint \psi_r(y,v) \nabla_\theta f(y, v) \rd y \rd v\,.
\end{align*}
Consequently $\sum_r \phi_r(x) \psi_r(y,v)$  enters the adjoint equation \eqref{eq:adjoint pde}, while  $H_\theta[ \nabla_\theta \df](x;\th)$ no longer appears in \eqref{grad}.

To obtain an efficient update formula
, in practice, we simplify the gradient by freezing the perturbation state and using the following approximation\footnote{We point out that the approximation \eqref{eq:approx grad} is constructed based on the small perturbation assumption $\df\ll 1$. However, it also turns out to work well even when the perturbation is moderately large. If the system perturbation grows rapidly, it is not always guaranteed that the convergence with the gradient approximation is stable, in which case one may resort to complete auto-differentiation. Nevertheless, our experiments also indicate that the convergence of the latter can be sometimes unstable either. Therefore, the appropriate optimization approach can be case-by-case.}:
\begin{equation}\label{eq:approx grad}
\nabla_\th \cJ(\th) \approx \int^T_0\iint (\nabla_{\th} H)[\df(t)](x;\th)\lambda(x,v,t)\partial_v f(x,v,t)\rd x \rd v \rd t.
\end{equation}
This approximation depends only on the adjoint variable and the explicit parameterization of the control kernel, enabling efficient updates using standard optimizers such as gradient descent or Adam. Our experiments indicate that the approximate gradient \eqref{eq:approx grad} performs well for moderate perturbations. Modern deep learning packages such as
\href{https://github.com/rtqichen/torchdiffeq}{torchdiffeq} may be employed to implement auto-differentiation based on semi-discrete ODEs.

To achieve high efficiency in long-term simulations, we will numerically solve the Vlasov-Poisson equations \eqref{eq:fb optim constraint} using the semi-Lagrangian method, which is fully explicit and not constrained by a CFL condition on the time step:
\begin{itemize}
    \item Compute $f^{(1)}(x,v) = f(x-\frac{\Dt}{2} v,v,t^n)$.
    \item Compute $E[f^{(1)}](x)$ by solving $\partial_x E[f^{(1)}] = \rho^{(1)}-\ion$.
    \item Compute $f^{(2)}(x,v) = f^{(1)}\big(x,v-(E[f^{(1)}]+H[\df^{(1)}])\Dt\big)$.
    \item Compute $f^{n+1}(x,v) = f^{(2)}(x-\frac{\Dt}{2} v,v)$.
\end{itemize}
The solution $f^{n+1}(x,v)$ at the final stage is then taken as the approximation to $f(x,v,t^{n+1})$. When computing the discrete solution on the mesh $\{(x_i,v_j)\}$,  we apply bilinear interpolation to obtain the nodal values on the departure grids, $\{(x_i-\frac{\Dt}{2}v_j,v_j)\}$ for the first and third stages and $\{\big(x_i,v_j-(E[f^{(1)}](x_i)+H[\df^{(1)}](x_i))\Dt\big)\}$ for the second stage. Similarly, the backward time-stepping of the adjoint equation \eqref{eq:adjoint pde} can be implemented with 
\begin{itemize}
    \item Compute $\lam^{(-1)}(x,v) = \lam(x+\frac{\Dt}{2} v,v,t^n)$.
    \item Compute $\lam^{(-2)}(x,v) = \lam^{(-1)}\big(x,v+\frac{\Dt}{2}(E[f^{n-\frac{1}{2}}]+H[\df^{n-\frac{1}{2}}])\big)$, $\displaystyle f^{n-\frac{1}{2}} = \frac{1}{2}(f^n+f^{n-1})$, $\displaystyle \df^{n-\frac{1}{2}} = \frac{1}{2}(\df^n+\df^{n-1})$ .
    \item Compute $\displaystyle \lam^{(-3)}(x,v) = \lam^{(-2)}(x,v)-\Dt \df^{n-\frac{1}{2}}(x,v)$.
     \item Compute $\lam^{(-4)}(x,v) = \lam^{(-3)}\big(x,v+\frac{\Dt}{2}(E[f^{n-\frac{1}{2}}]+H[\df^{n-\frac{1}{2}}])\big)$.
    \item $\lam^{n-1}(x,v) = \lam^{(-4)}(x+\frac{\Dt}{2} v,v)$.
\end{itemize}
In the following numerical tests, we employ the adjoint-based method described above, except for one bump-on-tail example where we use \texttt{torchdiff} to obtain time-independent control.

We emphasize that the adjoint-based method has practical advantages for the present PDE-constrained optimization problem. Note that direct automatic differentiation through the full time-stepping solver would require storing or recomputing a large, multi-time-level computational graph that involves interpolation, field solves, feedback evaluation, and repeated time-stepping. This can lead to substantial memory and computational costs when $N_t$ is large. In contrast, the adjoint formulation computes the gradient through one backward adjoint problem and a temporal integral involving the state, the adjoint state, and the derivative of the feedback operator. It therefore avoids backpropagation through the full time-stepping graph and makes the roles of the VP dynamics, the self-consistent field, and the feedback operator explicit.


\section{Numerical experiments}\label{sec:numerical exp}
\subsection{Computational setups}\label{sec:computational setups}
We examine the effectiveness of operator-based feedback control \eqref{eq:H LNO} through numerical computations. In all test cases, the computational domain is chosen as $(x,v)\in [0,10\pi]\times[-8,8]$, with periodic boundary conditions applied at $x = 0$ and $x = 10\pi$. The background ion density is set to $\ion = 1$, and the self-generated electric field is computed by
\[
\left\{
\begin{array}{l}
    E(x,t) = -\partial_x \Phi(x,t),  \\
    \\
    -\partial^2_x\Phi(x,t) = \rho(x,t)-1,  \ \rho(x,t) = \int f(x,v,t) \rd v\,,\\
    \\
    \Phi(0,t) = \Phi(10\pi,t) = 0.
\end{array}
\right.
\]
The discrete solutions are computed on the uniform mesh with $N_x = 100$ and $N_v = 200$ grids in the $x$ and $v$ directions, respectively.  The semi-Lagrangian time integration of the forward and backward equations \eqref{eq:fb optim constraint} and \eqref{eq:adjoint pde} are implemented with time step $\Dt = 0.2$.

To enforce exact periodic boundary conditions, rather than learning the basis, we fix the basis $\{\phi_k(\cdot)\}_{k=1}^{r}$ in \eqref{eq:H LNO} using $r=31$ trigonometric functions,
\[
\{\phi_k(x)\}^{31}_{k=1} = \left\{1,\sin(\tfrac{l}{5}x),\cos(\tfrac{l}{5}x)\right\}_{1\leq l \leq 15}.
\]
The integral kernel, $\psi^{NN}(x,v;\th_\psi)=(\psi_1(x,v),\psi_2(x,v),\cdots,\psi_{31}(x,v))$, remains a multi-output neural network mapping $\R^2$ to $\R^{31}$. Here we set up $\psi^{NN}$ using a 4-layer fully connected neural network with the distribution of neurons $2-64-32-31$. The output of each hidden layer is activated by the \texttt{ReLU} function before being fed to the next layer. In all experiments, the controller $H[\df]$ is optimized with respect to the running loss $\displaystyle \frac{1}{2}\int^T_0 ||f(t)-\feq||^2_2\rd t$ with $T = 30$, following the lines of Section \ref{sec:adjoint method}. 

To stabilize training in the early stages, we initialize the controller close to zero by setting the parameters of the output layer of $\psi^{NN}$ to zeros. The parameters $\th_\psi$ are then updated using the approximate gradient \eqref{eq:approx grad} over 3000 iterations. For efficient convergence, we first perform 200 steps of \texttt{Adagrad} as a preprocessing phase, followed by \texttt{Adam} iterations to accelerate training in the later stages. \footnote{We point out that the combination of Adagrad and Adam iterations is employed due to the fact that the Adagrad method, with fast gradient decay, converges quickly at the early stage but gets stuck afterwards. In contrast, the Adam method can be oscillatory under randomly initialized network parameters, but it helps accelerate iterations after Adagrad preprocessing. Thus, combining optimizers with different features would better ensure the stability and efficiency of training at different phases. The optimal choice of optimizers and/or combinations will be left for future discussions.}. 

To expand the exploration of the state space and prevent the controller from overfitting to data within the finite time horizon $t \in [0, T]$, we add a perturbation $f_p$ to the initial data $f_0$ at each iteration when computing the gradient,  defined as
\[
f_p(x,v) = \eps_p\sum^{k_p}_{k=1}\sum^{n_p}_{n=1}\omega_{kn}\hphi_k(x)\hh_n(v)\,.
\]
Here $\{\hphi(x)\}$ and $\{\hh_n(v)\}$ are orthonormal basis in $x$ and $v$ directions, respectively. The weight vector, $\omega$ = $\displaystyle \{\omega_{kn}\}_{1\leq k \leq k_p, 1\leq n\leq n_p}\in \R^{k_p\times n_p}$, is randomly and uniformly sampled from the $k_p\times n_p-$dimensional unit ball. In our computations, we construct $f_p$ using the normalized trigonometric functions on $x\in [0,10\pi]$,
\[
\{\hphi_k(x)\}^{11}_{k=1} = \left\{\sqrt{\frac{1}{10\pi}},\sqrt{\frac{1}{5\pi}}\sin(\tfrac{l}{5}x),\sqrt{\frac{1}{5\pi}}\cos(\tfrac{l}{5}x) \right\}_{1\leq l \leq 5},
\]
and the orthonormal Hermite basis functions on $v\in (-\infty,\infty)$,
\[
\hh_n(v) = \sqrt{\frac{1}{\sqrt{2\pi}n!}}He_n(v) \exp(-\tfrac{v^2}{4}), \quad 0\leq n \leq 5.
\]
Here, the Hermite polynomials $He_n(v)$, can be computed recursively,
\[
\begin{split}
&He_0(v) = 1,\\
&He_{n+1}(v) = vHe_n(v)-nHe_{n-1}(v), \quad n\geq 1.
\end{split}
\]

At the end of training, we use the future loss, $||f(2T)-\feq||_2$, as a measure of generalizability to select the best iteration. The parameters corresponding to this iteration are then adopted as the final trained model. This approach differs from simply taking the parameters at the end of the 3000 iterations. Because our training introduces noise, the iterations may explore multiple local minima rather than converge to a single one. Therefore, the training should be viewed as a search process, and the historically best result is used.

 To demonstrate the benefit of introducing dynamical feedback control, we will also include the time-independent counterpart for comparison, given by
\begin{equation} \label{Hindep}
H_{indep}(x) = \sum^{15}_{k=1}\theta_k\sin(\tfrac{k}{5}x)+\sum^{15}_{k=0}\theta_{k+16}\cos(\tfrac{k}{5}x).
\end{equation}
In this case, the gradient $\nabla_\theta H_{indep}(x;\theta)$ reduces to $(\phi_1(x),\phi_2(x),\cdots \phi_{31}(x))^\top$ in the gradient \eqref{eq:approx grad}. 

\begin{remark}
By adding a random initial perturbation $f_p$ at each iteration, the training procedure can be viewed as optimizing the expected trajectory loss over a distribution of perturbed initial data. Let $\mu$ be the distribution used to generate $f_p$. The corresponding stochastic optimization problem is
\[
\min_{\theta}\;
\mathbb E_{f_p\sim \mathcal \mu}
\bigl[l_T(\theta;f_p)\bigr],
\qquad
l_T(\theta;f_p)
=
\int_0^T
\|\delta f_\theta(t;f_p)\|_{L^2_{x,v}}^2\,dt .
\]
In principle, the expectation can be approximated by the empirical
average
\[
\frac1N\sum_{i=1}^N l_T(\theta;f_{p,i}).
\]
In the present implementation, the gradient at each optimization step is computed using one freshly sampled perturbation, i.e., $N=1$. We refer the interested reader to Appendix~\ref{append:stat err} for the associated statistical error analysis.
\end{remark}


\subsection{Two-stream instability}\label{sec:two stream}
We consider the following two-stream distribution as the target equilibrium:
\[
\feq(v) = \frac{1}{2\sqrt{2\pi}}\exp(
-\tfrac{(v-\bar{v})^2}{2}
)+\frac{1}{2\sqrt{2\pi}}\exp(
), \quad \bar{v} = 2.4.
\]
Assume that the Vlasov-Poisson system is initialized with the perturbed state
\[
f_0(x,v) = (1+\eps \cos(\tfrac{x}{5}))\feq(v), \quad x\in[0,10\pi], \quad \eps = 0.001.
\]
The two-stream equilibrium is known to be unstable, with perturbations amplified by the self-generated electric field. Our goal is to suppress this instability using an operator-based external field. The training procedure follows the approach outlined in Section \ref{sec:computational setups}, employing a 
learning rate of $lr = 5\times10^{-3}$ for the \texttt{Adagrad} iterations and $lr = 5\times10^{-4}$ for the subsequent \texttt{Adam} iterations. The magnitude of the random initial perturbation $f_p$ is set to $\eps_p=10^{-5}$.

\begin{figure}[h!]
    \centering
    \begin{subfigure}{0.45\textwidth}
    \includegraphics[width=\textwidth]{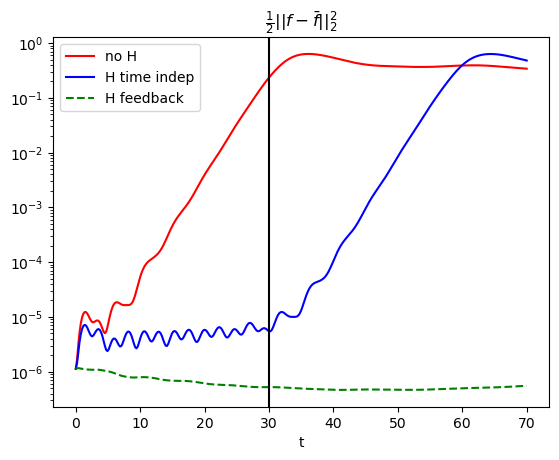}
    \subcaption{$L^2-$perturbation}
    \end{subfigure}
    \quad
    \begin{subfigure}{0.45\textwidth}
    \includegraphics[width=\textwidth]{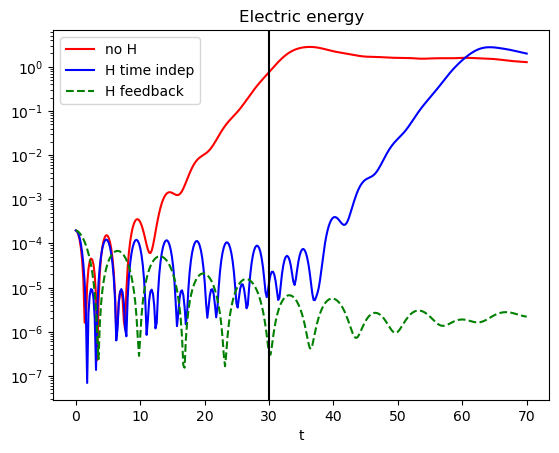}
    \subcaption{electric energy}
    \end{subfigure}
    \caption{Two-stream instability. History of $L^2-$state perturbation and electric energy over $t\in[0,70]$. No external field (red lines) vs time-independent control (blue lines) vs dynamical feedback control (green lines). The vertical black line marks the terminal time of optimization.}
    \label{fig:two stream perturb}
\end{figure}

\begin{figure}[h!]
    \centering
    \begin{subfigure}{0.3\textwidth}
    \includegraphics[width=\textwidth]{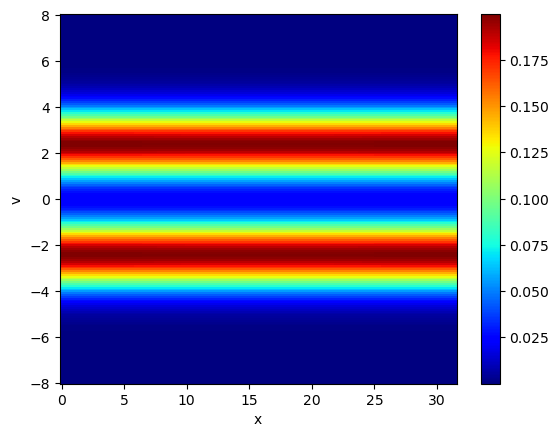}
    \subcaption{no $H$, $t=0$}
    \end{subfigure}
    \quad
    \begin{subfigure}{0.3\textwidth}
    \includegraphics[width=\textwidth]{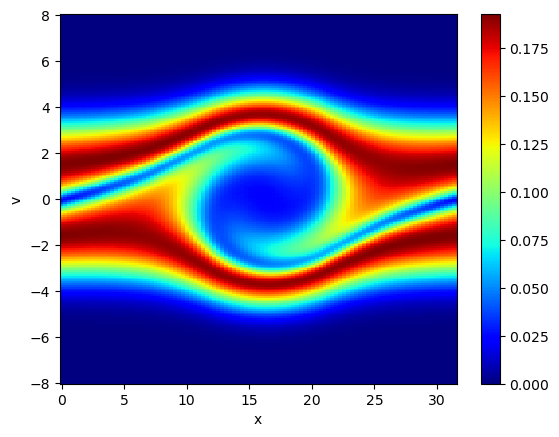}
    \subcaption{no $H$, $t=35$}
    \end{subfigure}
    \begin{subfigure}{0.3\textwidth}
    \includegraphics[width=\textwidth]{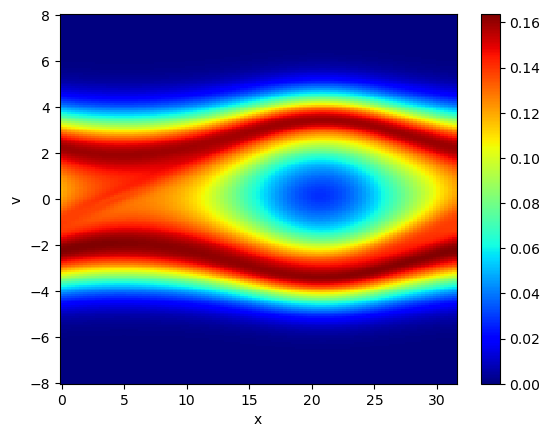}
    \subcaption{no $H$, $t=70$}
    \end{subfigure}
    \\
    \begin{subfigure}{0.3\textwidth}
    \includegraphics[width=\textwidth]{ts_f_t0.png}
    \subcaption{$H$ time-independent, $t=0$}
    \end{subfigure}
    \quad
    \begin{subfigure}{0.3\textwidth}
    \includegraphics[width=\textwidth]{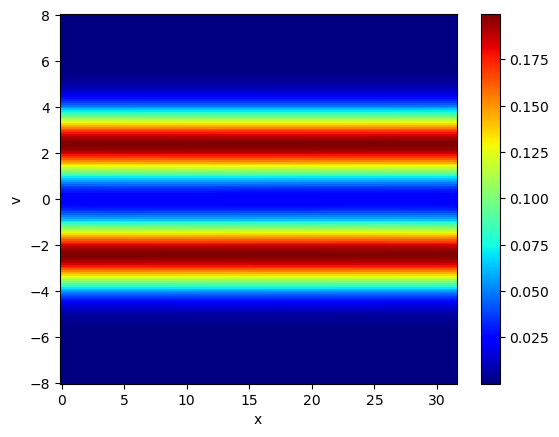}
    \subcaption{$H$ time-independent, $t=35$}
    \end{subfigure}
    \begin{subfigure}{0.3\textwidth}
    \includegraphics[width=\textwidth]{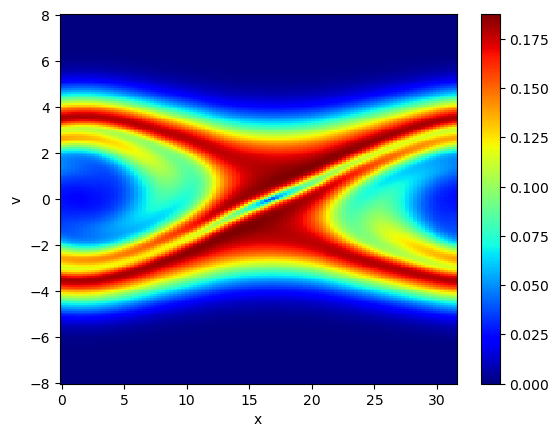}
    \subcaption{$H$ time-independent, $t=70$}
    \end{subfigure}
    \\
    \begin{subfigure}{0.3\textwidth}
    \includegraphics[width=\textwidth]{ts_f_t0.png}
    \subcaption{$H$ feedback, $t=0$}
    \end{subfigure}
    \quad
    \begin{subfigure}{0.3\textwidth}
    \includegraphics[width=\textwidth]{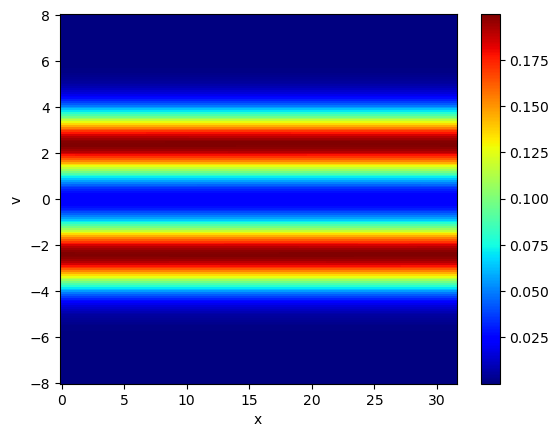}
    \subcaption{$H$ feedback, $t=35$}
    \end{subfigure}
    \begin{subfigure}{0.3\textwidth}
    \includegraphics[width=\textwidth]{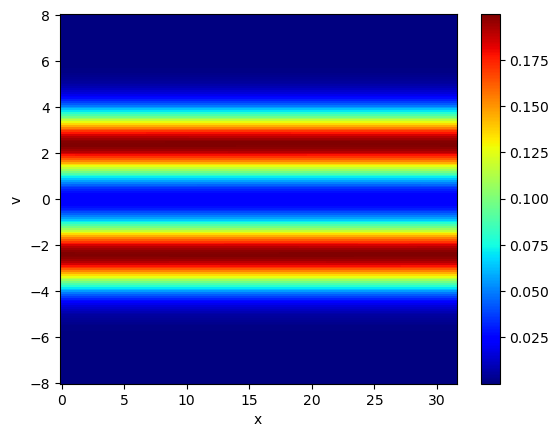}
    \subcaption{$H$ feedback, $t=70$}
    \end{subfigure}
    \caption{Two-stream instability. Evolution of $f(x,v,t)$. No external field (first row) vs time-independent control (second row) vs dynamical feedback control (third row).}
    \label{fig:two stream f}
\end{figure}

\begin{figure}[h!]
    \centering
    \includegraphics[scale = 0.45]{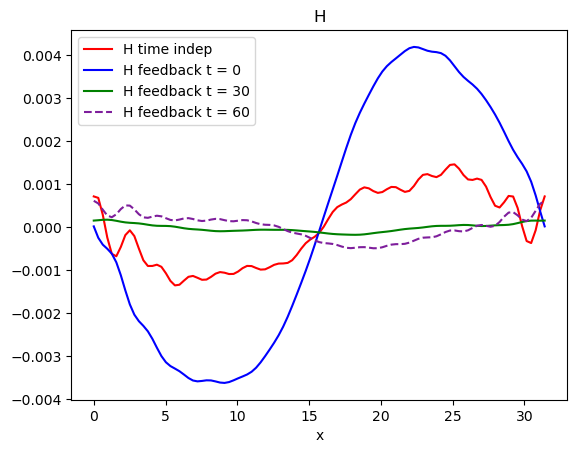}
    \caption{Two-stream instability. External fields by time-independent control and dynamical feedback control.}
    \label{fig:two stream H}
\end{figure}

Figure \ref{fig:two stream perturb} shows the evolution of the $L^2-$state perturbation, $\displaystyle \frac{1}{2}\iint |f(x,v,t)-\feq(v)|^2\rd x \rd v $, and the electric energy of the self-generated field, $\displaystyle\frac{1}{2}\int E(x,t)^2\rd x $, over the time interval $t\in[0,70]$. It is observed that the time-independent external field, optimized over $t\in[0,30]$, effectively suppresses the perturbation within the optimization interval but quickly loses effectiveness as time progresses. In contrast, the dynamical feedback control improves stability within the optimization interval and maintains effectiveness well beyond it. This behavior is further confirmed by the particle state distribution $f$, shown in Figure \ref{fig:two stream f}. 
Figure \ref{fig:two stream H} depicts the external fields at different time points, revealing that for $t>T$, the structure of the operator-based external field can be notably different from that of $H_{indep}$.

\subsection{Bump-on-tail instability}\label{sec:bump on tail}
\begin{figure}[h!]
    \centering
    \begin{subfigure}{0.45\textwidth}
    \includegraphics[width=\textwidth]{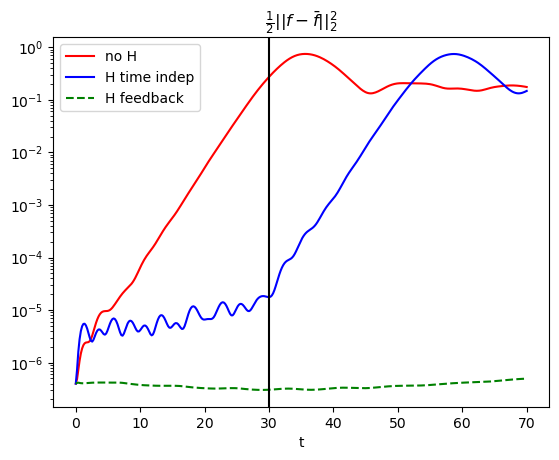}
    \subcaption{$L^2-$perturbation}
    \end{subfigure}
    \quad
    \begin{subfigure}{0.45\textwidth}
    \includegraphics[width=\textwidth]{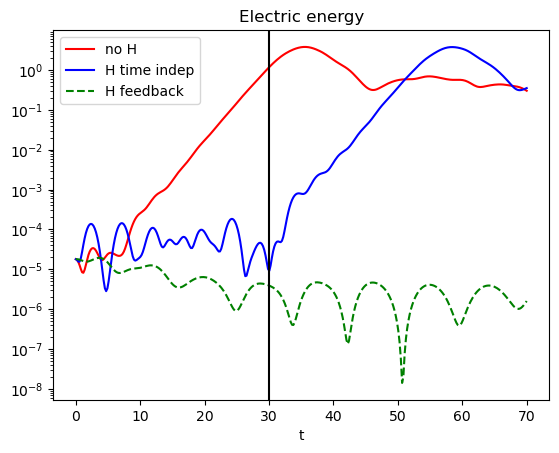}
    \subcaption{electric energy}
    \end{subfigure}
    \caption{Bump-on-tail instability. History of $L^2-$state perturbation and electric energy over $t=70$. No external field (red lines) vs time-independent control (blue lines) vs dynamical feedback control (green lines). The vertical black line marks the terminal time of optimization.}
    \label{fig:bump on tail perturb}
\end{figure}

\begin{figure}[h!]
    \centering
    \begin{subfigure}{0.3\textwidth}
    \includegraphics[width=\textwidth]{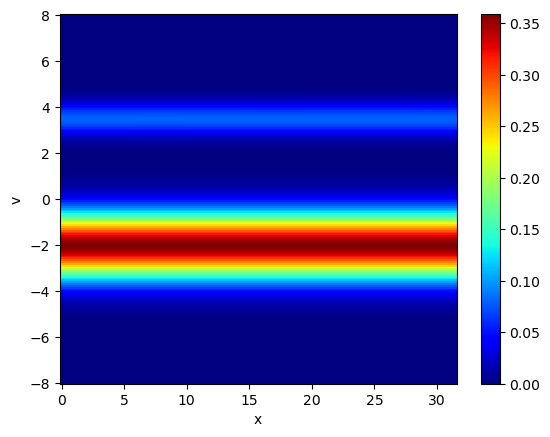}
    \subcaption{no $H$, $t=0$}
    \end{subfigure}
    \quad
    \begin{subfigure}{0.3\textwidth}
    \includegraphics[width=\textwidth]{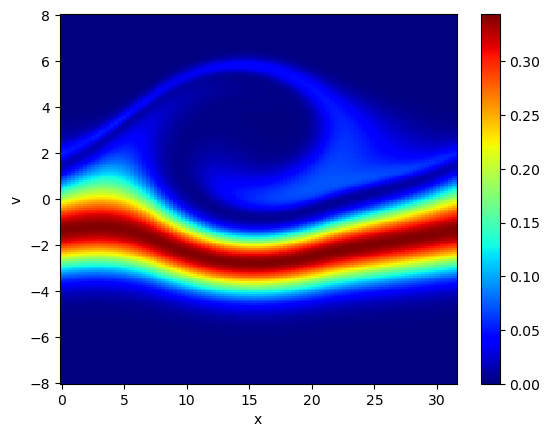}
    \subcaption{no $H$, $t=35$}
    \end{subfigure}
    \begin{subfigure}{0.3\textwidth}
    \includegraphics[width=\textwidth]{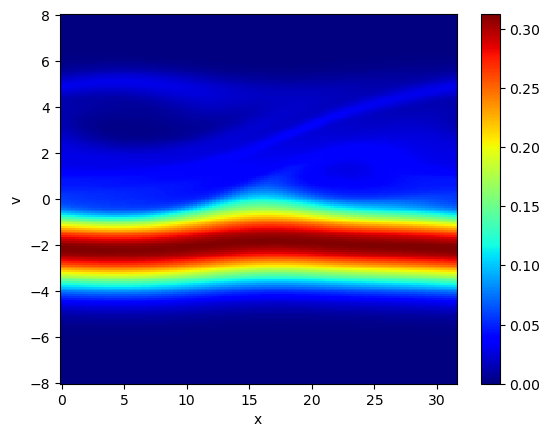}
    \subcaption{no $H$, $t=70$}
    \end{subfigure}
    \\
    \begin{subfigure}{0.3\textwidth}
    \includegraphics[width=\textwidth]{bt_f_t0_vt0.25.png}
    \subcaption{$H$ time-indep., $t=0$}
    \end{subfigure}
    \quad
    \begin{subfigure}{0.3\textwidth}
    \includegraphics[width=\textwidth]{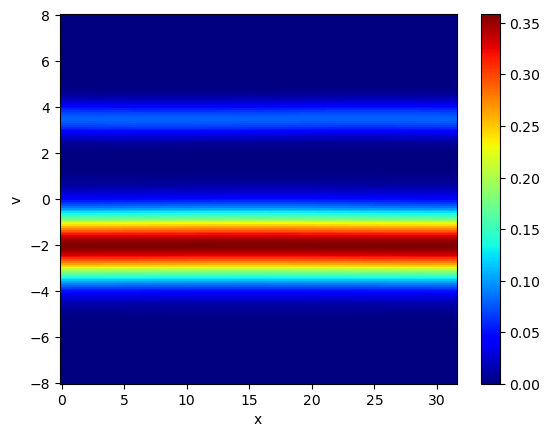}
    \subcaption{$H$ time-indep., $t=35$}
    \end{subfigure}
    \begin{subfigure}{0.3\textwidth}
    \includegraphics[width=\textwidth]{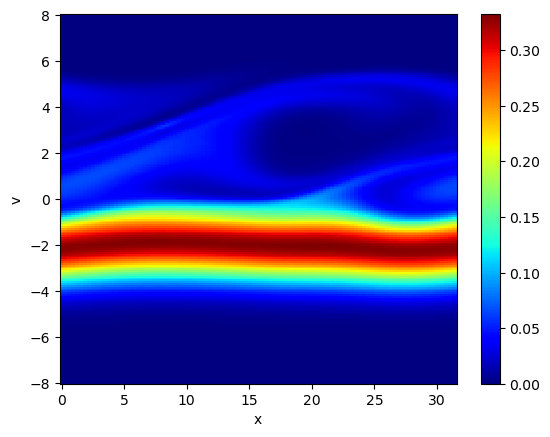}
    \subcaption{$H$ time-indep., $t=70$}
    \end{subfigure}
    \\
    \begin{subfigure}{0.3\textwidth}
    \includegraphics[width=\textwidth]{bt_f_t0_vt0.25.png}
    \subcaption{$H$ feedback, $t=0$}
    \end{subfigure}
    \quad
    \begin{subfigure}{0.3\textwidth}
    \includegraphics[width=\textwidth]{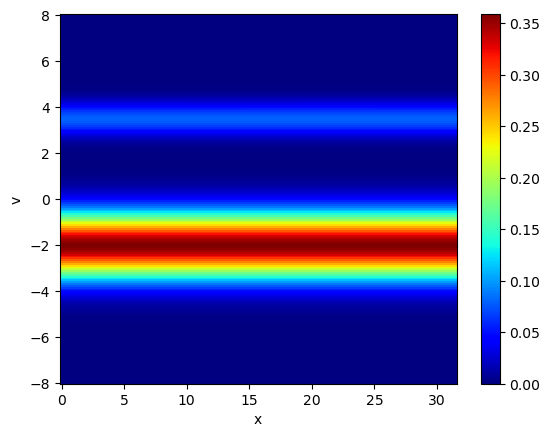}
    \subcaption{$H$ feedback, $t=35$}
    \end{subfigure}
    \begin{subfigure}{0.3\textwidth}
    \includegraphics[width=\textwidth]{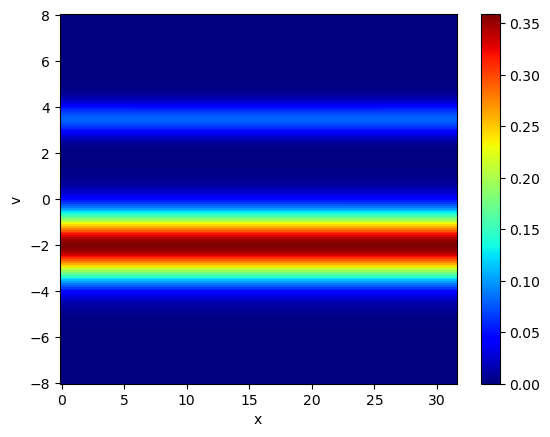}
    \subcaption{$H$ feedback, $t=70$}
    \end{subfigure}
    \caption{Bump-on-tail instability. Evolution of particle distribution $f(x,v,t)$. No external field (first row) vs time-independent control (second row) vs dynamical feedback control (third row).}
    \label{fig:bump on tail f}
\end{figure}

We also examine the performance of our operator-based feedback control when applied to suppress the bump-on-tail instability. The following target equilibrium is considered:
\[
\feq(v) = \frac{\omega_1}{\sqrt{2\pi}}\exp(-\tfrac{(v-\bar{v}_1)^2}{2})+\frac{\omega_2}{\sqrt{2\pi v_t}}\exp(-\tfrac{(v-\bar{v}_2)^2}{2v_t}).
\]
In our experiment, we set $\omega_1=0.9$, $\omega_2=0.1$, $\bar{v}_1 = -2$, $\bar{v}_2 = 3.5$, $v_t = 0.25$. We initialize the VP system by imposing a small perturbation to the equilibrium,
\[
f_0(x,v) = \feq(v)+\frac{\eps\omega_2}{\sqrt{2\pi v_t}}\exp(-\tfrac{(v-\bar{v}_2)^2}{2v_t})\sin(\tfrac{x}{5}), \quad x\in[0,10\pi], \quad \eps = 0.003,
\]
which leads to the destruction of the high-velocity thin tail as time develops. The controller network is trained with the learning rate $lr=2\times10^{-3}$ for \texttt{Adagrad} preprocessing and $lr=3\times10^{-4}$ for \texttt{Adam} iterations. 
The magnitude of the random initial perturbation $f_p$ is set to $\eps_p=10^{-6}$. 
We compare the solutions produced by the time-independent control \eqref{Hindep} and the feedback control \eqref{eq:H LNO} in Figures \ref{fig:bump on tail perturb}–\ref{fig:bump on tail f}, clearly demonstrating the superior effectiveness of the feedback control.

\begin{figure}[h!]
    \centering
    \includegraphics[scale = 0.45]{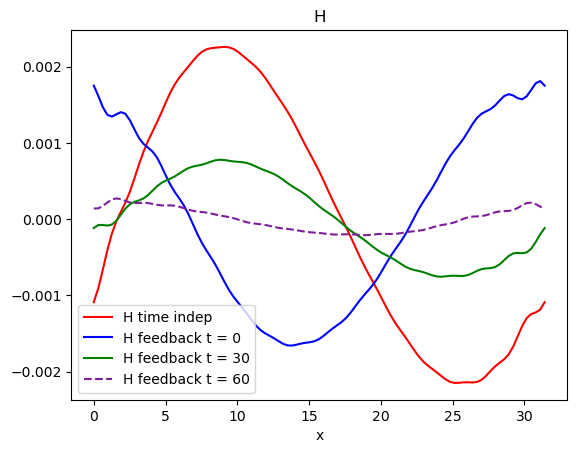}
    \caption{Bump-on-tail instability. External fields by time-independent control and dynamical feedback control.}
    \label{fig:bump on tail H}
\end{figure}

\subsection{Robustness under noisy feedback}\label{sec:noisy feedback}
In practical applications, measured data may be inaccurate, potentially undermining the effectiveness of feedback control. In this section, we test the robustness of our algorithm under noisy measurements. The operator-based controllers are the same as those obtained in Sections \ref{sec:two stream} and \ref{sec:bump on tail}, except that the external field, $H[\df^\sig(t)](x)$, perceives noisy feedback from the environment,
\begin{equation}\label{eq:noisy feedback}
\df^\sig(x,v,t) = \df(x,v,t)+\sig \mathcal{N}_{x,v,t}\,.
\end{equation}
Here $\{\mathcal{N}_{x,v,t}\}_{x,v,t}$ are independent and identically distributed according to the standard normal distribution.

\begin{figure}[h!]
    \centering
    \begin{subfigure}{0.3\textwidth}
    \includegraphics[width=\textwidth]{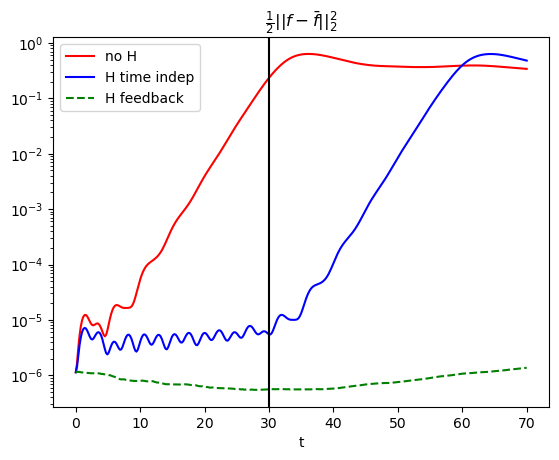}
    \subcaption{$\sig = 2\times 10^{-5}$ ($0.1\cdot ||\df_0||_\infty$)}
    \end{subfigure}
    \quad
    \begin{subfigure}{0.3\textwidth}
    \includegraphics[width=\textwidth]{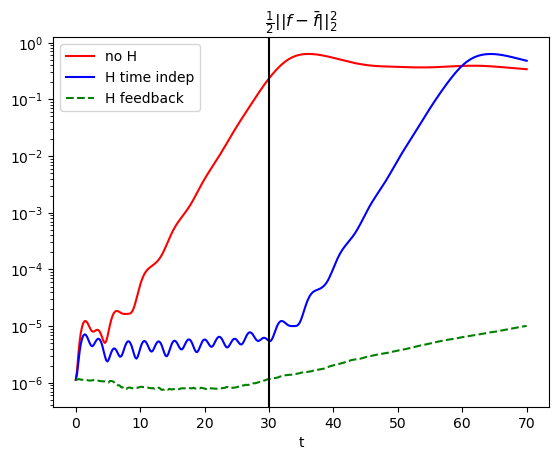}
    \subcaption{$\sig = 5\times 10^{-5}$ ($0.25\cdot||\df_0||_\infty$)}
    \end{subfigure}
    \quad
    \begin{subfigure}{0.3\textwidth}
    \includegraphics[width=\textwidth]{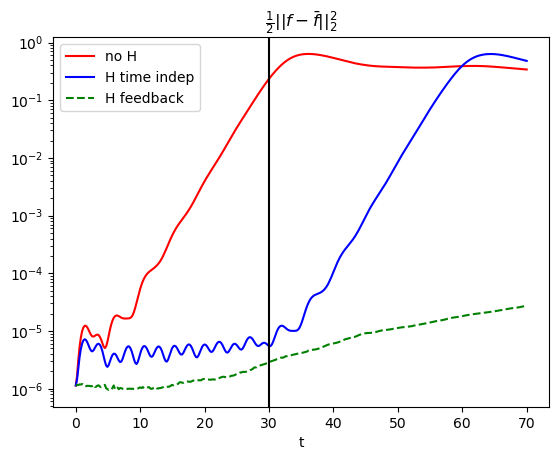}
    \subcaption{$\sig = 1\times 10^{-4}$ ($0.5\cdot||\df_0||_\infty$)}
    \end{subfigure}
    
    \caption{Two stream instability with noisy feedback. History of $L^2-$state perturbation over $t\in[0,70]$. No external field (red lines) vs time-independent control (blue line) vs dynamical feedback control (green lines) under varying noise magnitudes $\sig$. The vertical black line marks the terminal time of optimization.}
    \label{fig:two stream perturb noise}
\end{figure}

\begin{figure}[h!]
    \centering
    \begin{subfigure}{0.3\textwidth}
    \includegraphics[width=\textwidth]{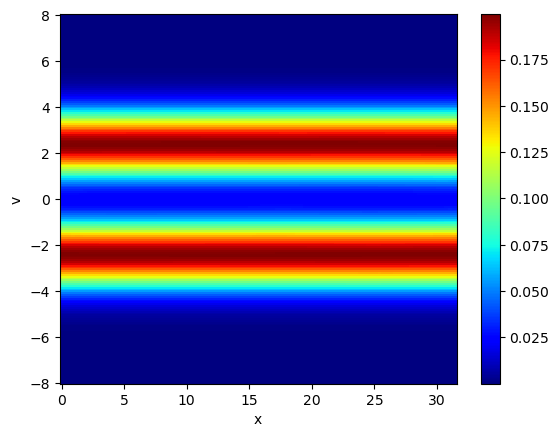}
    \subcaption{$\sig = 2\times 10^{-5}$ ($0.1\cdot ||\df_0||_\infty$)}
    \end{subfigure}
    \quad
    \begin{subfigure}{0.3\textwidth}
    \includegraphics[width=\textwidth]{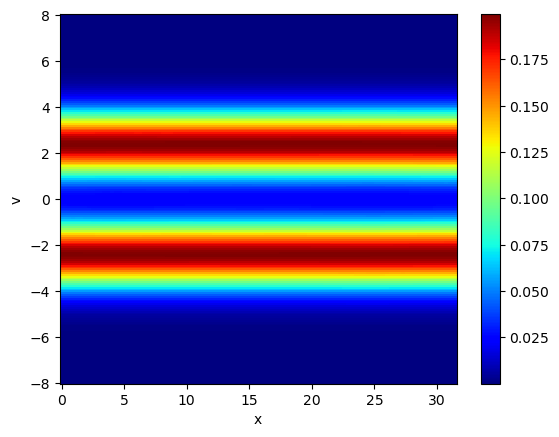}
    \subcaption{$\sig = 5\times 10^{-5}$ ($0.25\cdot ||\df_0||_\infty$)}
    \end{subfigure}
    \quad
    \begin{subfigure}{0.3\textwidth}
    \includegraphics[width=\textwidth]{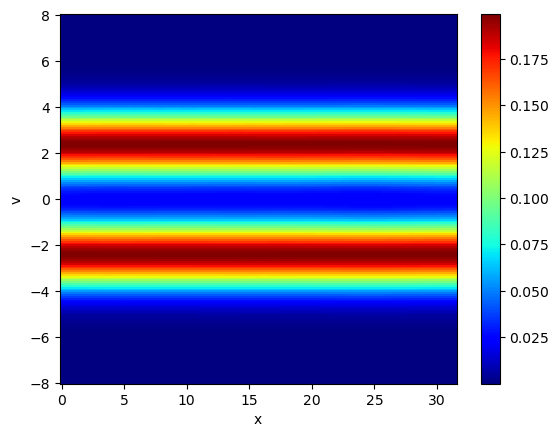}
    \subcaption{$\sig = 1\times 10^{-4}$ ($0.5 \cdot||\df_0||_\infty$)}
    \end{subfigure}
    
    \caption{Two stream instability with noisy feedback. Particle distribution $f(x,v,t)$ at $t=70$ generated by dynamical feedback control under varying noise magnitude $\sig$.}
    \label{fig:two stream f noise}
\end{figure}

\begin{figure}[h!]
    \centering
    \begin{subfigure}{0.3\textwidth}
    \includegraphics[width=\textwidth]{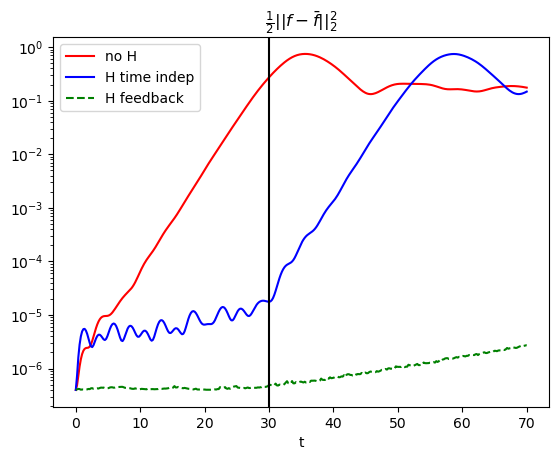}
    \subcaption{$\sig = 2.4\times 10^{-5}$ ($0.1\cdot ||\df_0||_\infty$)}
    \end{subfigure}
    \quad
    \begin{subfigure}{0.3\textwidth}
    \includegraphics[width=\textwidth]{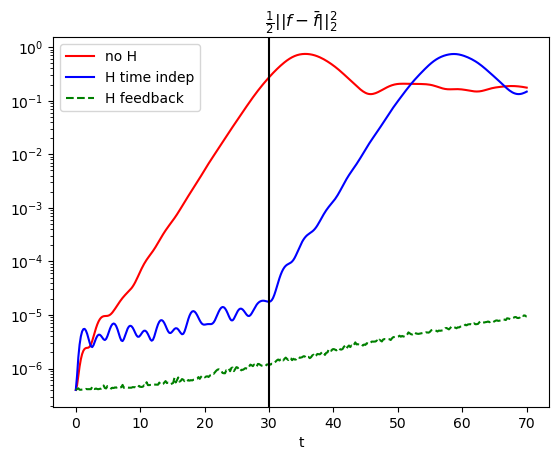}
    \subcaption{$\sig = 6\times 10^{-5}$ ($0.25\cdot||\df_0||_\infty$)}
    \end{subfigure}
    \quad
    \begin{subfigure}{0.3\textwidth}
    \includegraphics[width=\textwidth]{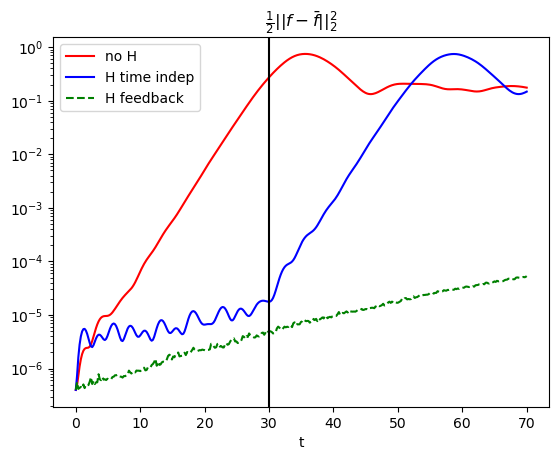}
    \subcaption{$\sig = 1.2\times 10^{-4}$ ($0.5\cdot ||\df_0||_\infty$)}
    \end{subfigure}
    
    \caption{Bump-on-tail instability with noisy feedback. History of $L^2-$state perturbation. No external field (red lines) vs time-independent control (blue line) vs dynamical feedback control (green lines) under varying noise magnitude $\sig$. The vertical black line marks the terminal time of optimization.}
    \label{fig:bump on tail perturb noise}
\end{figure}

\begin{figure}[h!]
    \centering
    \begin{subfigure}{0.3\textwidth}
    \includegraphics[width=\textwidth]{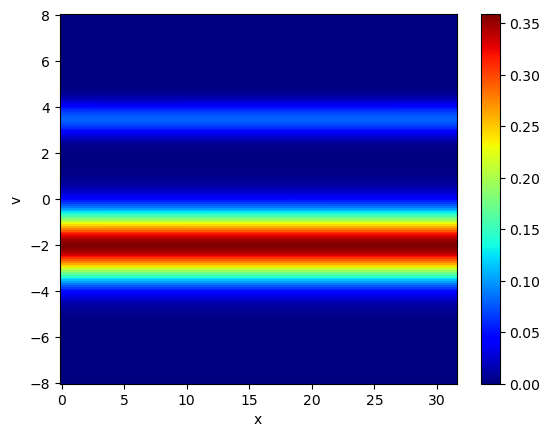}
    \subcaption{$\sig = 2.4\times 10^{-5}$ ($0.1\cdot ||\df_0||_\infty$)}
    \end{subfigure}
    \quad
    \begin{subfigure}{0.3\textwidth}
    \includegraphics[width=\textwidth]{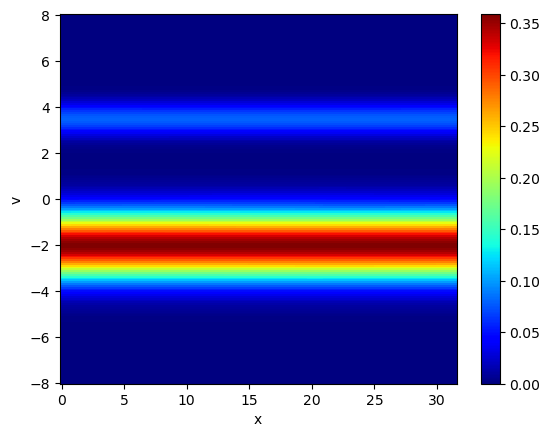}
    \subcaption{$\sig = 6\times 10^{-5}$ ($0.25\cdot||\df_0||_\infty$)}
    \end{subfigure}
    \quad
    \begin{subfigure}{0.3\textwidth}
    \includegraphics[width=\textwidth]{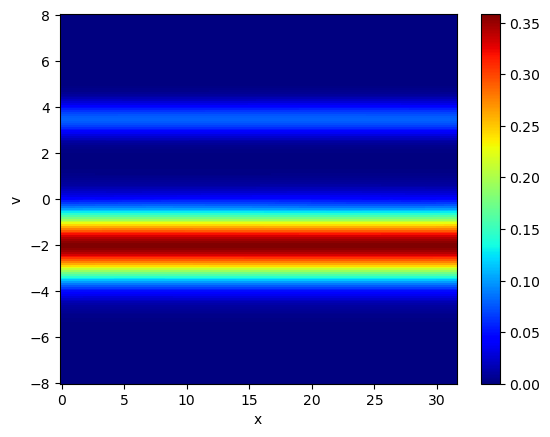}
    \subcaption{$\sig = 1.2\times 10^{-4}$ ($0.5\cdot ||\df_0||_\infty$)}
    \end{subfigure}
    
    \caption{Bump-on-tail instability with noisy feedback. Particle distribution $f(x,v,t)$ at $t=70$ generated by dynamical feedback control under varying noise magnitude $\sig$.}
    \label{fig:bump on tail f noise}
\end{figure}

Figures \ref{fig:two stream perturb noise} -- \ref{fig:two stream f noise} show simulations of the two-stream instability, with measurements $\df$ affected by noise of varying strengths. The initial particle state perturbation has magnitude $||\df_0||_\infty \approx 1.99\times 10^{-4}$. In Figure \ref{fig:two stream perturb noise} we observe that when the noise is small (compared to the initial perturbation), the operator-based external field maintains a clear advantage in terms of providing effective long-time instability control. As the magnitude of the noise increases to $\sig=1\times10^{-4}$ ($\approx 50\%  ||\df_0||_\infty$) the performance of $H[\df^\sig]$ is notably deteriorated. However, for $t>30$, our feedback control construction can still suppress the growth of the perturbation much better than the time-independent counterpart. As shown in Figure \ref{fig:two stream f noise}, the particle distribution generated by noisy feedback control remains relatively close to the target equilibrium up to $t=70$. Figures \ref{fig:bump on tail perturb noise} -- \ref{fig:bump on tail f noise} show the computations for the bump-on-tail instability ($||\df_0||_\infty \approx 2.39\times 10^{-4}$). Although the inaccurate measurement of perturbation weakens the performance of our algorithm, the control remains effective over a long time period, even under relatively large noise.

\begin{remark}
We point out that the magnitude $\eps_p$ of the artificial perturbations added to the initial condition during training, although small, can have a notable effect on the robustness of the learned feedback law under noisy observations. Figure~\ref{fig:two stream compare epsp} displays six independent successful training runs for each choice of $\eps_p$. Here, by a successful training run, we mean that, within the prescribed maximum number of optimization iterations, the learned controller achieves $||\delta f(2T)||_2<||\delta f_0||_2$ in the clean-feedback test.

Under clean observations, the two choices of $\eps_p$ yield comparable closed-loop performance, and all trained controllers effectively suppress the perturbation. Under stronger observation noise, however, the controllers trained with $\eps_p=10^{-5}$ exhibit more concentrated behavior: five of the six runs remain close to one another, while one converges to a controller with poorer noise generalization. In contrast, the noisy-feedback performance for $\eps_p=10^{-6}$ varies substantially across the six runs.

These experiments suggest that a larger training perturbation can improve the consistency of robustness to off-nominal observations and reduce the frequency of poorly generalizing controllers. Nevertheless, choosing $\eps_p$ too large may introduce greater variability into the training trajectories and make the PDE-constrained optimization more difficult to converge. In practice, $\eps_p$ should therefore be chosen to cover the expected perturbation and observation-noise regime while maintaining stable convergence of the training procedure.

\begin{figure}[h!]
    \centering
    \begin{subfigure}{0.45\textwidth}
        \includegraphics[width=\textwidth]{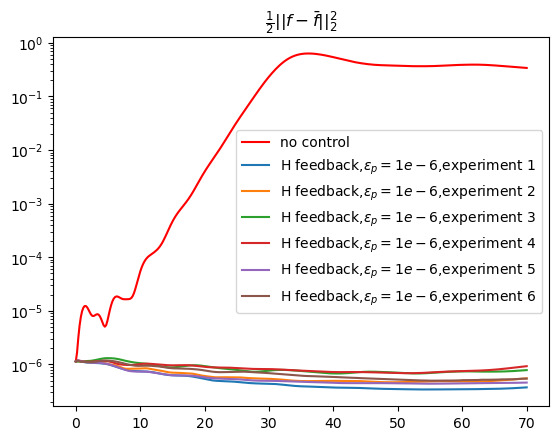}
        \subcaption{$\eps_p=10^{-6}$, $\sig=0$}
    \end{subfigure}
    \quad
    \begin{subfigure}{0.45\textwidth}
        \includegraphics[width=\textwidth]{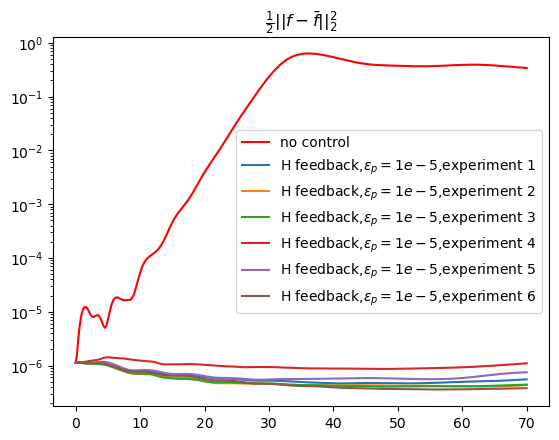}
        \subcaption{$\eps_p=10^{-5}$, $\sig=0$}
    \end{subfigure}
    \\
    \begin{subfigure}{0.45\textwidth}
        \includegraphics[width=\textwidth]{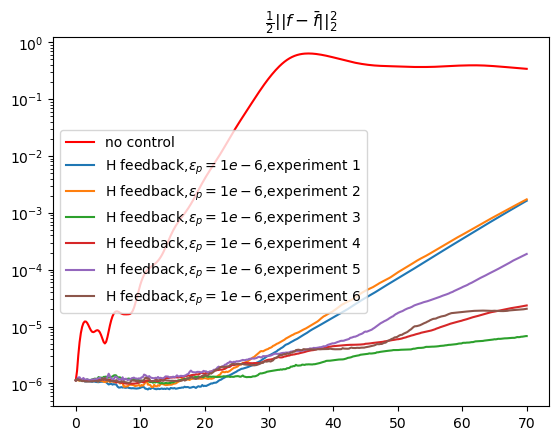}
        \subcaption{$\eps_p=10^{-6}$, $\sig=10^{-4}$}
    \end{subfigure}
    \quad
    \begin{subfigure}{0.45\textwidth}
        \includegraphics[width=\textwidth]{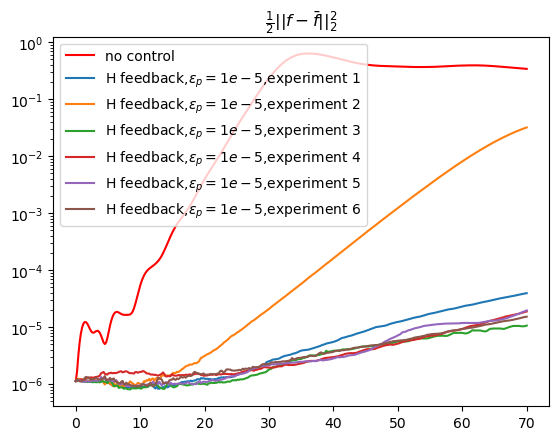}
        \subcaption{$\eps_p=10^{-5}$, $\sig=10^{-4}$}
    \end{subfigure}
    \caption{Two-stream instability: histories of the $L^2$ state perturbation for feedback controllers obtained from six independent converged training runs for each value of $\eps_p$. The upper row shows clean observations, while the lower row shows noisy observations with $\sig=10^{-4}$. Each controlled curve corresponds to one independently trained controller, and the red curve denotes the uncontrolled solution. The controllers trained with $\eps_p=10^{-5}$ exhibit more consistent robustness across the repeated runs under noisy observations.}
    \label{fig:two stream compare epsp}
\end{figure}
\end{remark}

\subsection{Effect of modal truncation}\label{sec:modal truncation}

 In practical settings, the externally imposed field may not be realizable as an arbitrary function of space, but only through a finite number of spatial modes or actuator channels. It is therefore useful to examine how many output modes are needed for the learned feedback law to remain effective.

We emphasize that the following test is not a retraining study over different neural operator ranks. Instead, we take the feedback controller trained with $K=15$ spatial Fourier modes (consistent with the tests in Sections \ref{sec:two stream}--\ref{sec:noisy feedback}) and, after the controller produces the field $H[\delta f](x)$, project this field onto the truncated Fourier space
\[
\mathrm{span}\{1,\sin(\frac{k}{5}x),\cos(\frac{k}{5}x):1\leq k\leq K_{\rm tr}\}.
\]
The truncated field is then used in the VP simulation. This isolates the effect of the number of realizable output control modes.

\begin{figure}[h!]
    \centering
    \begin{subfigure}{0.47\textwidth}
    \includegraphics[width=\textwidth]{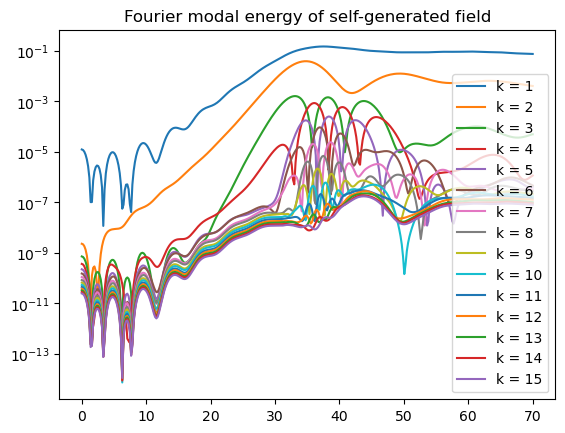}
    \subcaption{Modal energy of $E$}
    \end{subfigure}
    \quad 
    \begin{subfigure}{0.46\textwidth}
    \includegraphics[width=\textwidth]{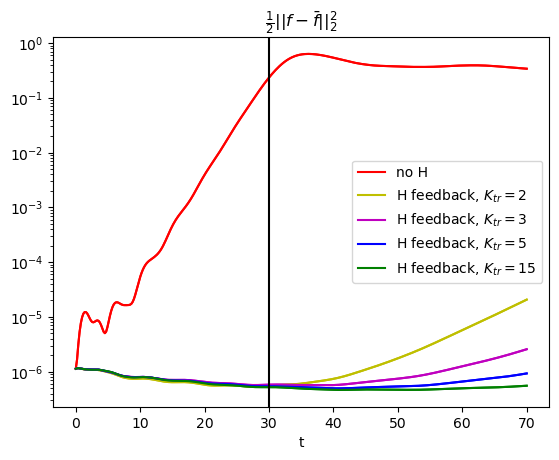}
    \subcaption{$L^2-$state perturbation}
    \end{subfigure}
    \caption{Modal truncation test for two-stream instability. The plot on the left shows the modal energy of the self-generated field $E$ in the absence of an external field. The plot on the right compares the growth of $L^2-$state perturbation when the controller output is truncated with different numbers of modes.  The vertical black line marks the terminal time of optimization.}
    \label{fig:two stream modal truncate}
\end{figure}

Figure~\ref{fig:two stream modal truncate} shows the modal energy of the self-consistent electric field in the uncontrolled two-stream dynamics. The growth is primarily concentrated in the first few spatial modes. In line with this observation, the feedback controller remains effective when the output field retains modes up to $K_{\rm tr}=5$, while more aggressive truncation to $K_{\rm tr}\leq 3$ leads to visible deterioration. This suggests that the retained control modes should cover the dominant unstable field modes, with a small number of additional modes to accommodate nonlinear mode interactions.

\begin{figure}[h!]
    \centering
    \begin{subfigure}{0.47\textwidth}
    \includegraphics[width=\textwidth]{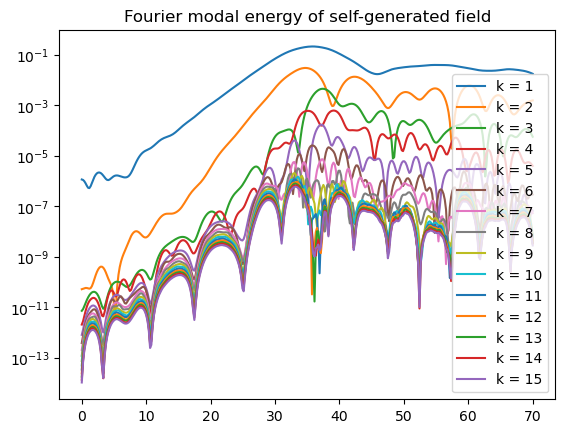}
    \subcaption{Modal energy of $E$}
    \end{subfigure}
    \quad 
    \begin{subfigure}{0.46\textwidth}
    \includegraphics[width=\textwidth]{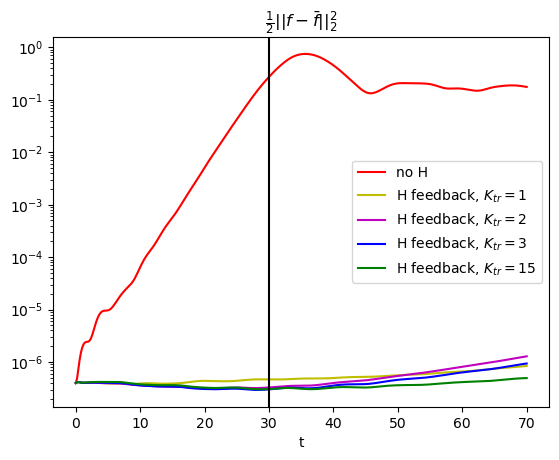}
    \subcaption{$L^2-$state perturbation}
    \end{subfigure}
    \caption{Modal truncation test for bump-on-tail instability. The plot on the left shows the modal energy of the self-generated field $E$ in the absence of an external field. The plot on the right compares the growth of $L^2-$state perturbation when the controller output is truncated with different numbers of modes.  The vertical black line marks the terminal time of optimization.}
    \label{fig:bump on tail modal truncate}
\end{figure}

The relation between the uncontrolled field spectrum and the required control modes, however, should be viewed only as a diagnostic guideline rather than a strict prescription. In the bump-on-tail example (Figure \ref{fig:bump on tail modal truncate}), the uncontrolled electric field develops a more distributed modal spectrum at later times, but the feedback controller remains effective even when the output field is truncated to the first spatial mode. This is consistent with the single-mode structure of the imposed perturbation in this test: suppressing the primary unstable channel early can sometimes prevent the subsequent nonlinear generation of higher modes. 

 In the main experiments, we use $K=15$ to provide a sufficiently rich modal control space without tuning the truncation level for each test case separately. But in practice, the required output dimension for sufficient control can be problem-dependent.

\section{A quasi-optimal universal control} \label{sec4}

\subsection{A cancellation-based construction}\label{sec:cancellation}

While the PDE-constrained optimization problem \eqref{eq:fb optim} offers a flexible framework for optimizing the controller of the Vlasov–Poisson system (and even more general plasma models such as the magnetically confined Vlasov–Maxwell equations), the optimal control for the nonlinear system \eqref{eq:vlasov poisson} can be sensitive to the choice of the initial dataset $f_0$ and the time horizon $T$.
 If the dataset is not exhaustive or $T$ is not sufficiently large, the resulting control may lose effectiveness over time, necessitating retraining and incurring significant computational cost.
To address this limitation, we propose a universal operator architecture that reduces such dependence. Specifically, we develop a quasi-optimal control strategy based on electric field cancellation, which depends only on the target equilibrium and can therefore be applied robustly across a broad range of initial conditions.

Subtracting the equilibrium $\feq(x,v)$, which satisfies $v\partial_x \feq+\baE\partial_v \feq = 0$, from the original VP system \eqref{eq:vlasov poisson}, we obtain the evolution equation for  $\df(x,v,t) = f(x,v,t)-\feq(x,v)$:
\begin{equation}\label{eq:vlasov poison perturb form}
    \partial_t \df+v\partial_x\df+\overline{E}\partial_v\df+(\dE+H)\partial_v f =0, 
\end{equation}
where $\baE(x)$ is the self-generated electric field of the equilibrium. The electric field perturbation, $\dE = E-\baE$, can be written in the operator form:
\[
\begin{split}
 \dE[\df(t)](x) =  -\partial_x(-\partial_x^2)^{-1}\drho(x,t),\quad \drho(x,t) = \int\df(x,v,t)\rd v .
\end{split}
\]
In particular, when the analytical form of the Green's function $G(x,y)$ is available for the Poisson equation with specific boundary condition, $\dE$ can be conveniently evaluated with
\begin{equation}\label{eq:dE operator form}
    \dE[\df(t)](x) = -\int\partial_x G(x,y)\drho(y,t)\rd y = -\iint\partial_x G(x,y)\df(y,v,t)\rd v \rd y .
\end{equation}
For instance, for the Poisson equation defined on the interval $x\in[a,b]$ with zero boundary data (which is the computational setup in Section \ref{sec:computational setups}), the Green's function is explicitly given by: 
\[
G(x,y) = \left\{\begin{array}{l}
     \displaystyle \frac{(x-a)(b-y)}{b-a}, \quad a\leq x \leq y \leq b \,, \\
     \\
     \displaystyle \frac{(y-a)(b-x)}{b-a}, \quad a\leq y \leq x \leq b\,.
\end{array}
\right.
\]
Alternatively, efficient algorithms such as fast Fourier transforms \cite[Chapter 8]{grossmann2007numerical}, multigrid methods \cite{zhang1998fast,guillet2011simple}, and fast multipole methods \cite{rokhlin1985rapid,ethridge2001new} can be applied to solve the Poisson equation directly.  

Physically, perturbations around an unstable equilibrium are amplified by the self-generated electric field. To counter this effect, we consider the construction
\[
H[\df(t)](x) = -\dE[\df(t)](x)+\dH[\df(t)](x).
\]
The goal is to cancel the destabilizing self-generated field $E$ while introducing an auxiliary term $\delta H$ to further suppress the perturbation. Substituting this external field into \eqref{eq:vlasov poison perturb form} and integrating the resulting equation against $\delta f$, we obtain
\begin{equation}\label{eq:df decay}
\begin{split}
\frac{1}{2}\frac{\rd}{\rd t}||\df(t)||^2_2 &= -\int\baE(x)\big[\int\df(x,v,t)\partial_v\df(x,v,t)\rd v \big] \rd x -\int \dH(x,t) \big[\int \df(x,v,t)\partial_v f(x,v,t)\rd v \big]\rd x \\
& = -\int \dH(x,t) \big[\int \df(x,v,t)\partial_v \feq(x,v)\rd v \big]\rd x .
\end{split}
\end{equation}
The second line is obtained by using the fact that $f = \feq+\df$ and $\int\df\partial_v\df \rd v = 0$. In particular, by setting 
\[
\dH[\df(t)](x) = \gam \int \df(x,v,t)\partial_v \feq(x,v)\rd v , \quad \gam>0,
\]
the equation \eqref{eq:df decay} implies the decay of perturbation in time,
\begin{align} \label{decay-df}
\frac{1}{2}\frac{d}{dt}||\df(t)||^2_2 = -\gam \int \big|\int \df(x,v,t)\partial_v \feq(x,v)\rd v \big|^2\rd x \leq 0.
\end{align}
This yields a convenient cancellation-based feedback controller,
\begin{equation}\label{eq:H cancellation based}
\begin{split}
& H[\df(t);\feq](x) = -\dE[\df(t)](x)+\gam \int \df(x,v,t)\partial_v \feq(x,v)\rd v , \quad \gam>0,\\
&\dE[\df(t)](x)  = -\iint\partial_x G(x,y)\df(y,v,t)\rd v \rd y .
\end{split}
\end{equation}
Note that this construction can be regarded as a special case of the linear integral operator \eqref{eq:H LNO} where both the basis and weights are fixed analytically. 
In general, however, exactly recovering this feedback law may require infinitely many spatial modes. When only finitely many control modes are retained, an additional finite-rank truncation error would be introduced. A standard Fourier-tail estimate for this error is provided in Appendix~\ref{append:finite rank err}. 

\begin{figure}[h!]
    \centering
    \begin{subfigure}{0.45\textwidth}
    \includegraphics[width=\textwidth]{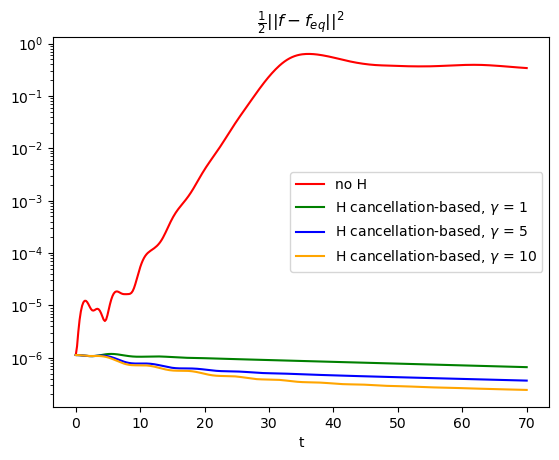}
    \subcaption{$L^2-$perturbation}
    \end{subfigure}
    \quad
    \begin{subfigure}{0.45\textwidth}
    \includegraphics[width=\textwidth]{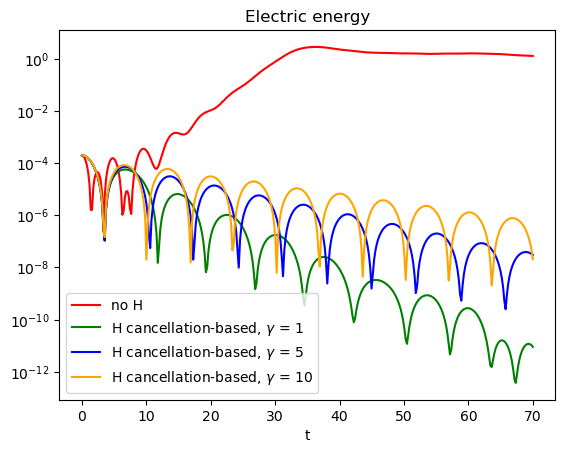}
    \subcaption{Electric field}
    \end{subfigure}
    \\
    \begin{subfigure}{0.45\textwidth}
    \includegraphics[width=\textwidth]{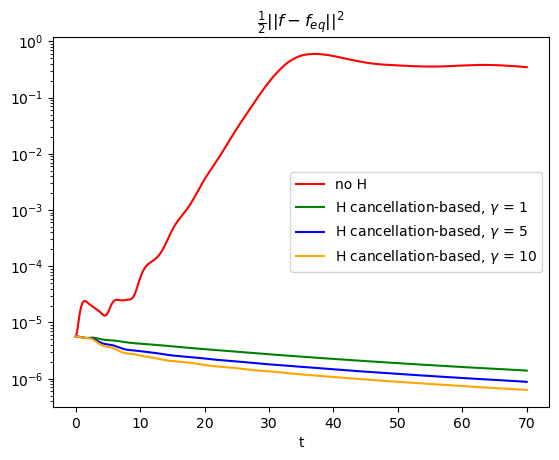}
    \subcaption{$L^2-$perturbation}
    \end{subfigure}
    \quad
    \begin{subfigure}{0.45\textwidth}
    \includegraphics[width=\textwidth]{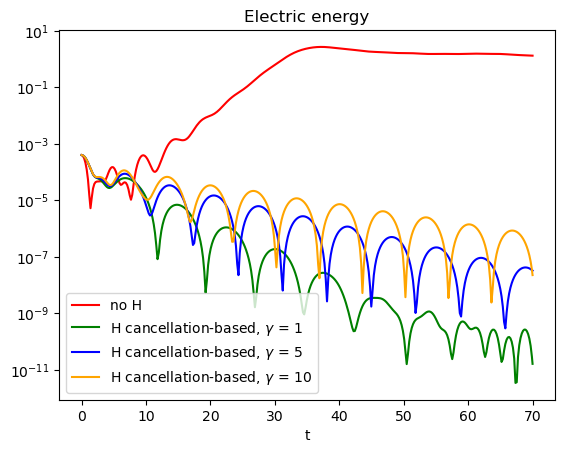}
    \subcaption{Electric field}
    \end{subfigure}
    \caption{Two-stream instability. Solutions generated by the cancellation-based controller \eqref{eq:H cancellation based} under different initial data. Solutions in the first row are generated with $f_0 = (1+0.001\cos(\frac{1}{5}x))\feq(v)$. Solutions in the second row are generated with $f_0 = (1-0.001\sin(\frac{1}{5}x)+0.002\cos(\frac{2}{5}x))\feq(v)$. }
    \label{fig:two stream cancellation}
\end{figure}

\begin{figure}[h!]
    \centering
    \begin{subfigure}{0.4\textwidth}
    \includegraphics[width=\textwidth]{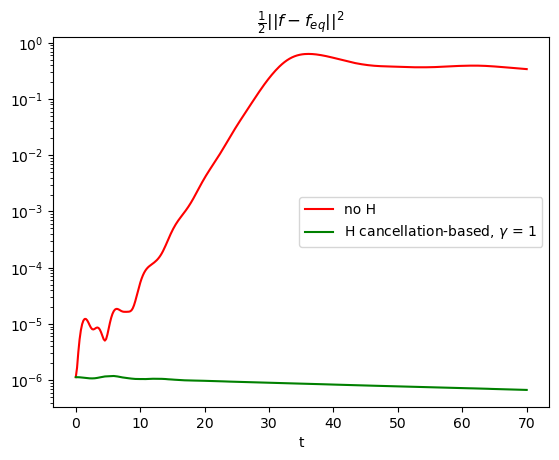}
    \subcaption{$\sig = 2\times 10^{-5}$ ($0.1\cdot ||\df_0||_\infty$)}
    \end{subfigure}
    \quad
    \begin{subfigure}{0.4\textwidth}
    \includegraphics[width=\textwidth]{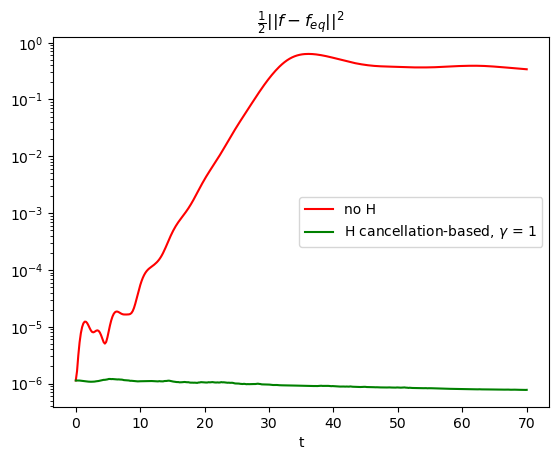}
    \subcaption{$\sig = 5\times 10^{-5}$ ($0.25\cdot||\df_0||_\infty$)}
    \end{subfigure}
    \\
    \begin{subfigure}{0.4\textwidth}
    \includegraphics[width=\textwidth]{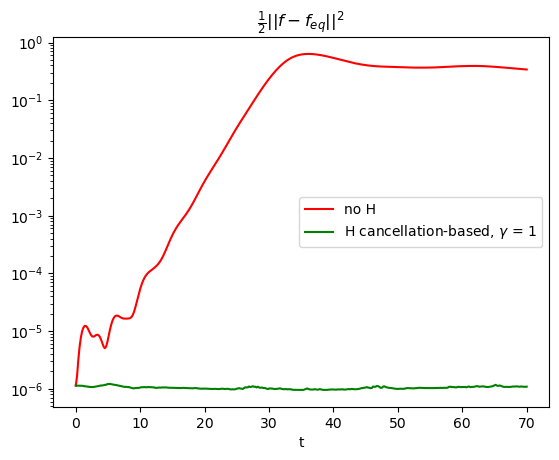}
    \subcaption{$\sig = 1\times 10^{-4}$ ($0.5\cdot||\df_0||_\infty$)}
    \end{subfigure}
    \quad
    \begin{subfigure}{0.4\textwidth}
    \includegraphics[width=\textwidth]{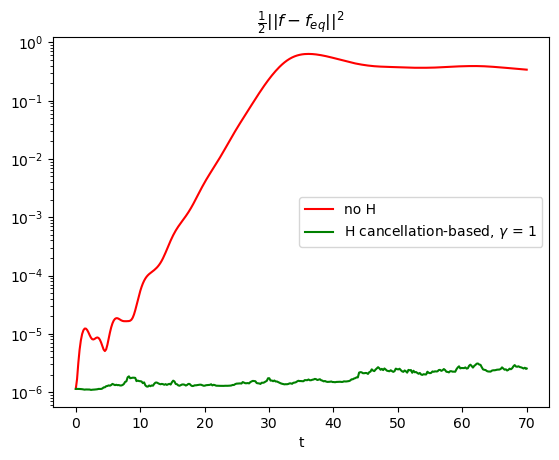}
    \subcaption{$\sig = 2\times 10^{-4}$ ($1\cdot||\df_0||_\infty$)}
    \end{subfigure}
    \caption{Two stream instability with noisy feedback. History of $L^2-$state perturbation over $t\in[0,70]$. No external field (red lines) vs cancellation-based control with $\gam = 1$ (green lines) under varying noise magnitudes $\sig$.}
    \label{fig:two stream cancellation noise}
\end{figure}

Although the control \eqref{eq:H cancellation based} does not ensure $||\df(t)||\overset{t\rightarrow \infty}{\longrightarrow}0$, it still provides satisfactory control when the initial perturbation is already small enough.
Figure \ref{fig:two stream cancellation} presents the evolution of the $L^2-$state perturbation $\frac{1}{2}||\df(t)||^2_2$, and the electric energy $\frac{1}{2}||E(t)||^2_2$, generated by the controller \eqref{eq:H cancellation based}  in the two-stream instability test case (see Section \ref{sec:two stream} for setups). The results confirm the universality of the controller with respect to different initial states $f_0$. Meanwhile, a larger $\gam$ leads to faster decay in the perturbation. Figure \ref{fig:two stream cancellation noise} demonstrates the performance of the controller under noisy feedback of the form \eqref{eq:noisy feedback}. The initial data is set to $f_0 = (1+0.001\cos(\frac{1}{5}x))\feq(v)$, in consistency with the setting in Section \ref{sec:noisy feedback}. The results show the notably improved robustness of the cancellation-based method as compared to the aforementioned neural operator-based construction (see Figure \ref{fig:two stream perturb noise}).

\begin{remark}
The construction in \eqref{eq:H cancellation based} is not unique. For example, an alternative choice is: 
\begin{equation}\label{eq:H cancellation based ver2}
\dH(x,t) = \gam \frac{\int|\df(x,v,t)|^2\rd v }{\int \df(x,v,t)\partial_v\feq(x,v)\rd v +\eps},
\end{equation}
where a small bias term $\eps$, having the same sign as $\int \df(x,v,t)\partial_v\feq(x,v)\rd v $, is introduced to mitigate severe round-off errors when the denominator is close to zero. Theoretically, the construction \eqref{eq:H cancellation based ver2} yields a sharper decay estimate,
\[
\frac{1}{2}\frac{d}{dt}||\df(t)||^2_2 = -\gam ||\df(t)||^2_2 \leadsto ||\df(t)||_2\leq ||\df(0)||^2_2 \cdot e^{-\gam t}.
\]
However, our numerical experiments show that the controller defined in \eqref{eq:H cancellation based ver2}, owing to its more intricate structure, is more sensitive to noisy feedback than the simpler linear controller in \eqref{eq:H cancellation based} 
(see Figure \ref{fig:two stream linear vs nonlinear})
. Consequently, the latter is preferable in practice, as it offers lower computational cost and greater robustness in noisy environments. 

\begin{figure}[h!]
    \centering
    \begin{subfigure}{0.4\textwidth}
    \includegraphics[width=\textwidth]{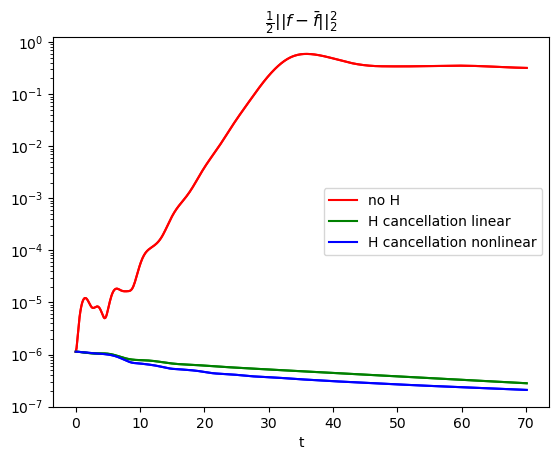}
    \subcaption{clean data ($\sig = 0$)}
    \end{subfigure}
    \quad 
    \begin{subfigure}{0.4\textwidth}
    \includegraphics[width=\textwidth]{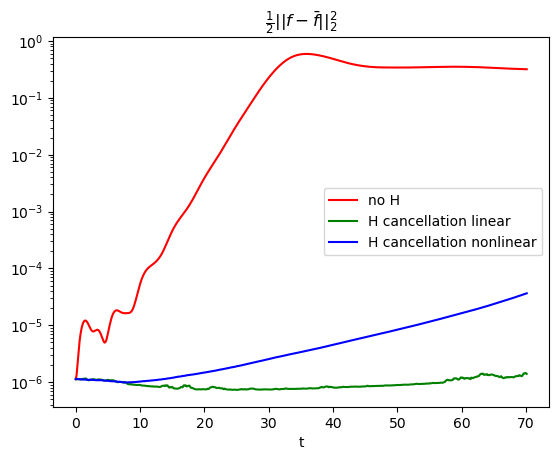}
    \subcaption{noisy feedback ($\sig = 1.5\times 10^{-4}$)}
    \end{subfigure}
    \caption{Two-stream instability. History of $L^2-$state perturbation generated with no control (red lines), linear control \eqref{eq:H cancellation based} with $\gam = 2$ (green lines), and nonlinear control \eqref{eq:H cancellation based ver2} with $\gam = 2$ (blue lines). Solutions are generated with $f_0 = (1+0.001\cos(\frac{1}{5}x))\feq(v)$ and different noise magnitudes $\sig$. Note that the nonlinear control generates slightly stronger decay in the clean data setting but is less robust under noisy feedback.}
    \label{fig:two stream linear vs nonlinear}
\end{figure}
\end{remark}

\begin{remark}
We highlight the difference between our cancellation-based feedback law and the analytical control field obtained in \cite{einkemmer2024control}, which was derived using a pole-elimination technique. This technique relies on the linearized system, along with a Fourier transform-based method. In contrast, our approach is formulated directly on the original nonlinear system and ensures analytical decay of $||\df(t)||_2$,  a property that is not guaranteed in \cite{einkemmer2024control}, particularly when the dynamics extend beyond the linear regime.
\end{remark}

\subsection{Extension to higher dimensions}

It is worth noting that the derivation \eqref{eq:vlasov poison perturb form}–\eqref{eq:H cancellation based} does not rely on any specific dimensionality. Hence, the operator construction \eqref{eq:H cancellation based} can be directly extended to higher-dimensional phase spaces, namely,
\begin{equation}\label{eq:H cancellation based multi-d}
\begin{split}
H[\df(t)](\bx) &= -\dE[\df(t)](\bx)+\gam\mathbb{P}_\nabla\big(\int_{\R^d} \df(t)\nabla_\bv \feq d\bv\big)(\bx),\\
\dE[\df(t)](\bx) &= \nabla \Delta^{-1}\drho(\bx,t).
\end{split}
\end{equation}
Here, the Helmholtz projection, $\mathbb{P}_\nabla(\bu) = \nabla\Delta^{-1}(\nabla\cdot\bu)$, removes the divergence-free portion of a vector field, ensuring that the dissipation term is an electric (potential) field\footnote{Note that $\int \df(\bx,\bv,t)\nabla_\bv\feq(\bx,\bv)d\bv$ is not an electric (potential) field in general. Thus, we perform Helmholtz decomposition $\int \df(\bx,\bv,t)\nabla_\bv\feq(\bx,\bv) d\bv = \nabla\phi_H(\bx,t)+\mathbf{c}_H(\bx,t)$ with 
\[\left\{\begin{array}{l}
    \Delta \phi_H =  \nabla \cdot (\int \df\nabla_\bv\feq d\bv), \quad \bx\in\Omega\\
    \\
     \phi_H = 0, \quad \bx\in \partial \Omega
\end{array}\right. \quad \text{and} \quad \left\{\begin{array}{l}
    \nabla\cdot\mathbf{c}_H = 0, \quad \bx\in\Omega\\
    \\
     \mathbf{c}_H\cdot\bn = (\int \df\nabla_\bv\feq d\bv)\cdot\bn, \quad \bx\in\partial\Omega
\end{array}\right..
\]
Then we take the potential part only to construct the dissipation field, $\dH = \nabla\phi_H$.}. The perturbation decay estimate is modified to 
\[
\frac{1}{2}\frac{d}{dt}||\delta f(t)||^2_2 = -\gam \int |\mathbb{P}_\nabla(\int \delta f(t)\nabla_\bv\feq d\bv)(x)|^2dx.
\]
Compared to the one-dimensional situation, the added complexity arises from the numerical evaluation of $\delta E$, velocity integrals, and the projection of the dissipation field. This universality makes the cancellation-based feedback strategy particularly appealing for plasma control problems beyond the one-dimensional setting.

\begin{figure}[h!]
\centering
    \begin{subfigure}{0.4\textwidth}
    \includegraphics[width=\textwidth]{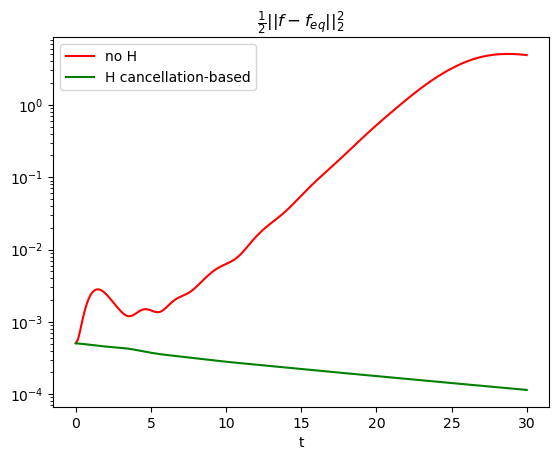}
    \subcaption{$L^2-$perturbation}
    \end{subfigure}
    \quad
    \begin{subfigure}{0.4\textwidth}
    \includegraphics[width=\textwidth]{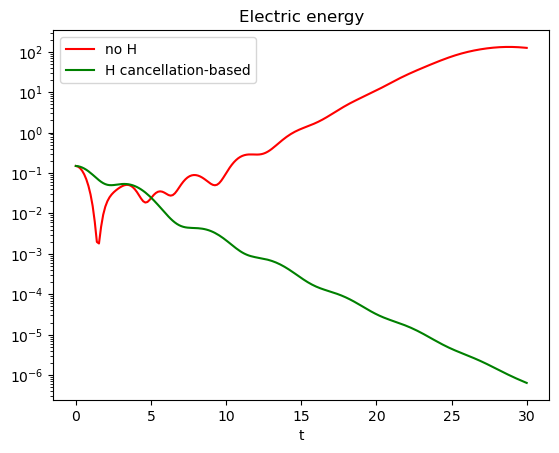}
    \subcaption{Electric energy}
    \end{subfigure}
    \caption{2-D two-stream instability. History of $L^2-$state perturbation, $\frac{1}{2}||\df(t)||^2_2$, and electric energy, $\int |E(\bx,t)|^2d\bx$, over $t\in[0,30]$. No external field (red lines) vs cancellation-based control \eqref{eq:H cancellation based multi-d} with $\gam = 2$ (green lines).}
    \label{fig:vp4d results}
\end{figure}

\begin{figure}[h!]
    \centering
    \begin{subfigure}{0.31\textwidth}
    \includegraphics[width=\textwidth]{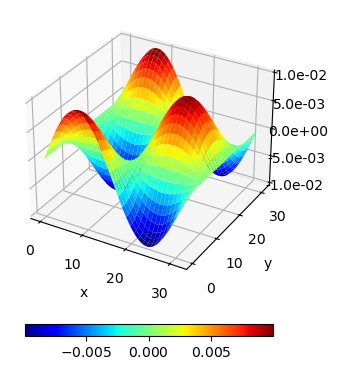}
    \subcaption{$t = 0$}
    \end{subfigure}
    \quad
    \begin{subfigure}{0.31\textwidth}
    \includegraphics[width=\textwidth]{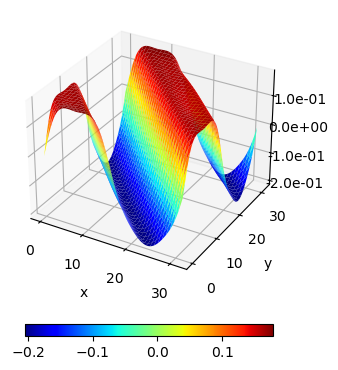}
    \subcaption{$t = 30$, no $H$}
    \end{subfigure}
    \quad
    \begin{subfigure}{0.31\textwidth}
    \includegraphics[width=\textwidth]{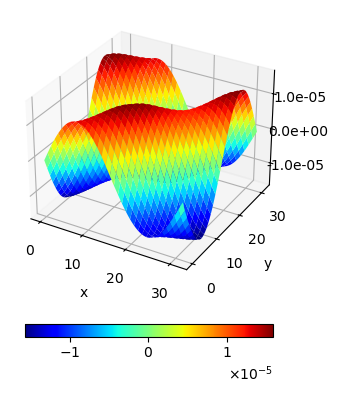}
    \subcaption{$t = 30$, $H$ cancellation}
    \end{subfigure}
    \\
    \begin{subfigure}{0.31\textwidth}
    \includegraphics[width=\textwidth]{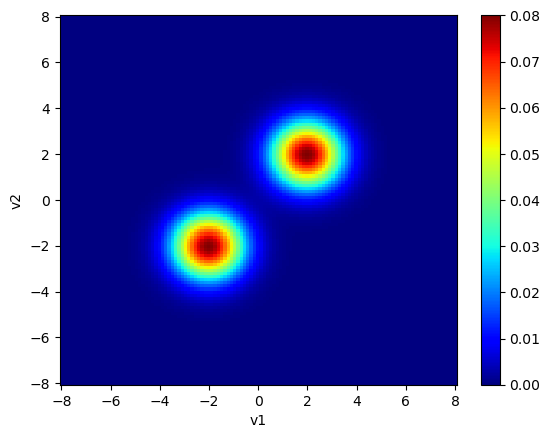}
    \subcaption{$t = 0$}
    \end{subfigure}
    \quad
    \begin{subfigure}{0.31\textwidth}
    \includegraphics[width=\textwidth]{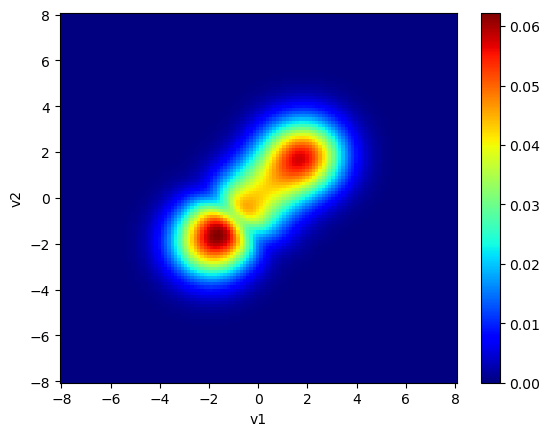}
    \subcaption{$t = 30$, no $H$}
    \end{subfigure}
    \quad
    \begin{subfigure}{0.31\textwidth}
    \includegraphics[width=\textwidth]{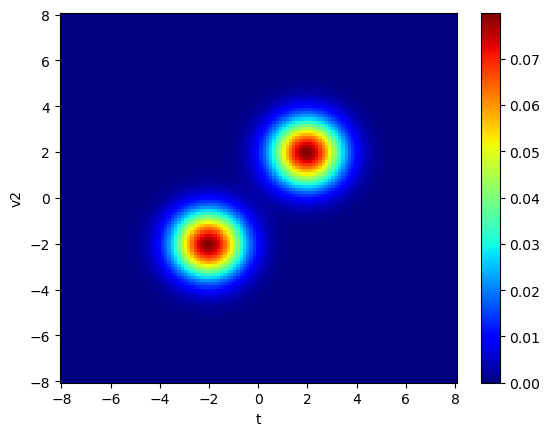}
    \subcaption{$t = 30$, $H$ cancellation}
    \end{subfigure}
    \caption{2-D two-stream instability. The first row displays surface plots of the density perturbation $\drho$. The different ranges of the color bars should be noted. The second row displays the contours of velocity distribution $f(6\pi,5\pi,v_1,v_2,t)$. The solution of the cancellation-based control \eqref{eq:H cancellation based multi-d} is obtained with $\gam = 2$.}
    \label{fig:vp4d contours}
\end{figure}

In Figures \ref{fig:vp4d results} -- \ref{fig:vp4d contours}, we present the simulation for the two-stream instability in two dimensions, where the two-dimensional Vlasov-Poisson equations,
\begin{equation}
 \left\{
\begin{array}{l}
\partial_t f(\bx,\bv,t)+\bv\cdot\nabla_\bx f(\bx,\bv,t)+(E(\bx,t)+H(\bx,t))\nabla_\bv f(\bx,\bv,t) = 0\,, \\
\\
E(\bx,t) = -\nabla \Phi(\bx,t), \quad \nabla\cdot E(\bx,t) = -\Delta \Phi(\bx,t) = \rho(\bx,t)-1 \,,\\
\\
\rho(\bx,t) = \int_{\R^2} f(\bx,\bv,t) d\bv\,,
\end{array}\right.
\end{equation}
are solved in the domain $\bx = (x,y)\in [0,10\pi]^2$, $\bv = (v_1,v_2)\in [-8,8]^2$. The equations are discretized using the semi-Lagrangian method with mesh grids $N_x = N_y = 70$, $N_{v_1} = N_{v_2} = 120$ and time step $\Dt = 0.15$. Periodic boundary conditions are applied at the four boundaries of the $xy-$region. The target equilibrium is set to 
\[
\feq(\bv) = \frac{1}{4\pi}\exp(-\tfrac{|\bv-\overline{\bv}|^2}{2})+\frac{1}{4\pi}\exp(-\tfrac{|\bv+\overline{\bv}|^2}{2}), \quad \overline{\bv} = (2,2),
\]
and the system is initialized with the perturbed data,
\[
f_0(\bx,\bv) = \big(1+\eps\sin(\tfrac{x}{5})\cos(\tfrac{y}{5}) \big)\feq(\bv), \quad \eps = 0.01.
\]
The results demonstrate that in the absence of control, the perturbation grows rapidly, whereas under the control \eqref{eq:H cancellation based multi-d}, the system is stabilized effectively in the 2D setting.

\begin{remark}
We emphasize that the derivation of \eqref{eq:H cancellation based} (or \eqref{eq:H cancellation based multi-d}) does \emph{not} rely on any assumption about 
$\feq$	or on linearization of the system. The only requirement is access to $\nabla_\bv \feq(\bx,\bv)$ which can be precomputed and stored. Although our experiments focus on a spatially invariant equilibrium $\feq(\bv)$to minimize numerical artifacts in perturbation growth, the method is, in principle, applicable to general anisotropic equilibria. Furthermore, the cancellation-based control can be regarded as a specific instance of the linear operator structure \eqref{eq:H LNO}, whereas the latter—with learnable basis functions and integral kernels—offers substantially greater flexibility. Consequently, we expect the linear LNO-based controller to generalize naturally to spatially inhomogeneous settings, provided the operator is trained with sufficient accuracy.
\end{remark}

\section{Conclusion}
In this paper, we develop a dynamic feedback control strategy for the Vlasov–Poisson system to address the challenge of long-time instability control in fusion energy applications. The core idea is to construct an operator that maps state perturbations to an external control field. We propose two approaches for constructing such an operator.

The first approach involves learning the operator using a neural network. Rather than relying on an off-the-shelf neural operator, which would impose a significant training burden, we draw inspiration from the linearized system and use optimal control theory to uncover the underlying structure of the ideal operator. This insight guides the design of a low-rank neural operator architecture. To train this model, we derive an adjoint-based gradient computation method. Compared to automatic differentiation, this adjoint approach is often more efficient and yields more accurate gradients, reducing the overall training cost. In the second approach, by explicitly analyzing the evolution equation satisfied by perturbations, we perform a direct energy estimate and propose a cancellation-based control strategy that eliminates the destabilizing component of the electric field. This leads to a novel closed-form operator, completely removing the need for training. The resulting operator is highly robust, even under noisy feedback.

Several interesting directions lie ahead. One concerns the theoretical understanding of the cancellation-based control defined in \eqref{eq:H cancellation based}. While the decay estimate \eqref{decay-df} suggests that the perturbation cannot be completely eliminated, since components orthogonal to $\partial_v \feq$	
are not directly controlled, numerical experiments nonetheless demonstrate effective overall decay. This raises the intriguing question of whether there is an underlying structure yet to be uncovered that explains the effectiveness of this operator. In particular, it would be interesting to show that if the initial perturbation is orthogonal to $\partial_v \feq$, then this orthogonal component still diminishes over time under the evolution governed by \eqref{eq:vlasov poison perturb form}.
Another promising direction involves the learning-based approach. To ensure long-term effectiveness of the learned operator, it is crucial to train it on a sufficiently rich set of perturbations, especially those likely to arise dynamically over time. This motivates the development of strategies for self-generating representative perturbations during training. Finally, an even more intriguing avenue is to extend the proposed strategy to the design of external magnetic fields. Since magnetic fields are often easier to tune in practice, this extension could significantly enhance the feasibility of real-world control implementations.


\section*{Acknowledgment}
JL and JC were supported by NSF-CCF:2212318, and JC was additionally supported by an Albert and Dorothy Marden Professorship. 
LW is supported in part by NSF grant DMS-1846854, DMS-2513336, and the Simons Fellowship. 

\appendix

\section{Statistical error induced by stochastic gradient}\label{append:stat err}

In this section, we provide a detailed analysis of the optimization error due to the empirical approximation of the stochastic gradient. Let $z = f_p$ denote an initial perturbation sampled from a probability distribution $\mu$ over an admissible perturbation set. For a feedback operator $H_\theta$, define the finite-horizon trajectory loss
\begin{equation}
\ell_T(\theta;z)
=
\int_0^T \|\delta f_\theta(t;z)\|_{L^2_{x,v}}^2\,\rd t,
\end{equation}
where $\delta f_\theta(t;z)$ is the solution of the closed-loop Vlasov--Poisson system initialized at $z$. The population and empirical risks are
\begin{equation}
\mathcal R_T(\theta)
=
\mathbb E_{z\sim\mu}\,\ell_T(\theta;z),
\qquad
\widehat{\mathcal R}_{N,T}(\theta)
=
\frac1N\sum_{i=1}^N \ell_T(\theta;z_i),
\end{equation}
where $z_1,\ldots,z_N$ are independent samples from $\mu$.

\begin{proposition}
Assume that the parameter satisfies $\theta \in \Theta\subset B_R(0)\subset\mathbb R^p$, and that for all $\theta,\theta'\in\Theta$ and all admissible $z$,
\begin{equation} \label{assump0}
0\le \ell_T(\theta;z)\le B_T,
\qquad
|\ell_T(\theta;z)-\ell_T(\theta';z)|
\le L_T\|\theta-\theta'\|.
\end{equation}
Then, for any $\varepsilon\in(0,R]$, with probability at least $1-\delta$,
\begin{equation*}
\sup_{\theta\in\Theta}
\left|
\mathcal R_T(\theta)
-
\widehat{\mathcal R}_{N,T}(\theta)
\right|
\le
2L_T\varepsilon
+
B_T
\sqrt{
\frac{
p\log(3R/\varepsilon)+\log(2/\delta)
}{2N}
}.
\label{eq:sampling-error-bound}
\end{equation*}
In particular, choosing $\varepsilon = \frac{R}{N}$, the sampling error scales as
\begin{equation*}
\sup_{\theta\in\Theta}
\left|
\mathcal R_T(\theta)
-
\widehat{\mathcal R}_{N,T}(\theta)
\right|
=
O\!\left(
B_T
\sqrt{
\frac{p\log N+\log(1/\delta)}{N}
}
\right),
\end{equation*}
up to constants depending on $R$, $L_T$, and $B_T$.
\end{proposition}
\begin{proof}
Let $\Theta_\varepsilon$ be an $\varepsilon$-net of $\Theta$. Since $\Theta\subset B_R(0)\subset\mathbb R^p$, the covering number satisfies
\begin{equation*}
|\Theta_\varepsilon|
\le
\mathcal N(\Theta,\varepsilon)
\le
\left(\frac{3R}{\varepsilon}\right)^p.
\end{equation*}
Here $\mathcal N(\Theta,\varepsilon)$ is the smallest number of radius-$\varepsilon$ balls needed to cover $\Theta$. 
For each fixed $\theta_j\in\Theta_\varepsilon$, the random variables
\begin{equation*}
X_i=\ell_T(\theta_j;z_i)
\end{equation*}
are independent and bounded in $[0,B_T]$. Hoeffding's inequality gives, for every $u>0$,
\begin{equation*}
\mathbb P\left(
\left|
\mathcal R_T(\theta_j)
-
\widehat{\mathcal R}_{N,T}(\theta_j)
\right|>u
\right)
\le
2\exp\left(-\frac{2Nu^2}{B_T^2}\right).
\end{equation*}
Taking a union bound over all $\theta_j\in\Theta_\varepsilon$ yields
\begin{equation*}
\mathbb P\left(
\max_{\theta_j\in\Theta_\varepsilon}
\left|
\mathcal R_T(\theta_j)
-
\widehat{\mathcal R}_{N,T}(\theta_j)
\right|>u
\right)
\le
2|\Theta_\varepsilon|
\exp\left(-\frac{2Nu^2}{B_T^2}\right).
\end{equation*}
Choosing
\begin{equation*}
u
=
B_T
\sqrt{
\frac{\log|\Theta_\varepsilon|+\log(2/\delta)}{2N}
}
\end{equation*}
shows that, with probability at least $1-\delta$,
\begin{equation*}
\max_{\theta_j\in\Theta_\varepsilon}
\left|
\mathcal R_T(\theta_j)
-
\widehat{\mathcal R}_{N,T}(\theta_j)
\right|
\le
B_T
\sqrt{
\frac{\log|\Theta_\varepsilon|+\log(2/\delta)}{2N}
}.
\label{eq:net-concentration}
\end{equation*}

Now take any $\theta\in\Theta$. Since $\Theta_\varepsilon$ is an $\varepsilon$-net, there exists $\theta_j\in\Theta_\varepsilon$ such that
\begin{equation*}
\|\theta-\theta_j\|\le \varepsilon.
\end{equation*}
By the Lipschitz assumption,
\begin{equation*}
|\mathcal R_T(\theta)-\mathcal R_T(\theta_j)|
\le
L_T\varepsilon,
\qquad
|\widehat{\mathcal R}_{N,T}(\theta)
-
\widehat{\mathcal R}_{N,T}(\theta_j)|
\le
L_T\varepsilon.
\end{equation*}
Therefore,
\begin{align*}
\left|
\mathcal R_T(\theta)
-
\widehat{\mathcal R}_{N,T}(\theta)
\right|
&\le
|\mathcal R_T(\theta)-\mathcal R_T(\theta_j)|
+
|\mathcal R_T(\theta_j)-\widehat{\mathcal R}_{N,T}(\theta_j)|
\notag \\
&\quad
+
|\widehat{\mathcal R}_{N,T}(\theta_j)
-
\widehat{\mathcal R}_{N,T}(\theta)|
\notag \\
&\le
2L_T\varepsilon
+
B_T
\sqrt{
\frac{\log|\Theta_\varepsilon|+\log(2/\delta)}{2N}
}.
\end{align*}
Taking the supremum over $\theta\in\Theta$ and using
\begin{equation*}
\log|\Theta_\varepsilon|
\le
p\log(3R/\varepsilon)
\end{equation*}
gives the result.
\end{proof}

\begin{remark}
    We comment on the assumptions used above. The assumption that $\Theta\subset B_R(0)$ is imposed in order to apply the $\varepsilon$-net argument and to simplify the estimate of the covering number. This assumption can be relaxed by working instead with Rademacher complexity. The assumption \eqref{assump0} may be justified through a stability estimate for the solution of the Vlasov-Poisson system with respect to perturbations in the external field, for which the result in \cite{loeper2006uniqueness} could be useful. We leave the detailed analysis to future work.
\end{remark}

\section{Modal truncation error for cancellation-based control}\label{append:finite rank err}

As discussed in Section~\ref{sec:cancellation}, the cancellation-based
control
\[
H[\delta f(t);\feq](x) = -\delta E[\delta f(t)](x)+\gamma\int \delta f(x,v,t)\partial_v \feq(v)\,dv, \qquad \gamma>0,
\]
yields a provable decay estimate for the perturbation energy. However, an exact representation of this feedback law in the operator form \eqref{eq:H LNO} may require infinitely many spatial modes. To quantify the effect of finite-mode truncation, we restrict ourselves to the one-dimensional case with $x\in[0, L]$ and fix the output basis functions as
\begin{align}\label{basis}
\phi_k(x) = \frac{1}{\sqrt L}e^{i\kappa_k x}, \qquad \kappa_k=\frac{2\pi k}{L}.
\end{align}
The integral kernels $\{\psi_k(y,v;\theta_\psi)\}$ are then the only learnable components.

The following proposition estimates how well a finite-rank operator of the form \eqref{eq:H LNO} can approximate the cancellation-based control.

\begin{proposition}\label{prop:cancellation finite rank}
Let $H[\delta f]$ denote the cancellation-based control
\eqref{eq:H cancellation based}, and let $H_r[\delta f]$ be the best finite-mode approximation of $H[\delta f]$ within the operator class \eqref{eq:H LNO}, where the output basis functions
$\{\phi_k(x;\theta_\phi)\}_{|k|\leq r}$ are fixed as in \eqref{basis}. Then we have the following estimate
\begin{align}\label{eq:finite rank err}
\|H_r [\delta f] - H [\delta f] \|_{L_x^2} \leq |\kappa_{r+1}|^{-1}\inner{\kappa_{r+1}}^{-s}
\|\drho\|_{H_x^s} + \gamma\langle \kappa_{r+1}\rangle^{-s}
\|\partial_v\bar f\|_{L_v^2}\,\|\delta f\|_{H_x^sL_v^2}\,,
\end{align}
where $\langle \kappa \rangle = \sqrt{1+ |\kappa|^2}$.
\end{proposition}

\begin{proof}
By standard approximation theory, given an orthonormal basis $\{\phi_k(x)\}_{|k|\le r}$, the best approximation of $H(x)$ in the $L^2$ sense over the space $\operatorname{span}\{\phi_k(x)\}_{|k|\leq r}$ is its orthogonal projection onto this subspace. Namely,
\begin{align}
H_r[\delta f](x) := P_r H[\delta f](x) = \sum_{|k|\leq r} \phi_k(x) \hat H_k, \quad \text{where} \quad  \hat H_k = \int \phi_k(x) H[\delta f](x) \rd x\,.
\end{align} 
Moreover, we have 
\begin{align} \label{eq:rank truncation err}
\|H- H_r \|_{L_x^2}^2 = \sum_{|k|>r}|\hat H_k|^2 \le
\langle \kappa_{r+1}\rangle^{-2s}
\sum_{|k|>r}\langle \kappa_k\rangle^{2s}|\hat H_k|^2 \leq 
\langle \kappa_{r+1}\rangle^{-2s}\|H\|_{H_x^s}^2.
\end{align}
Note that this estimate holds for a general $L^2_x$ function.

Consider $H[\delta f]$ in the form \eqref{eq:H cancellation based}. It has a Poisson part that leads to $\delta E$, and an additional damping part. We estimate them individually. For the Poisson part, write 
\[
\drho(x)=\sum_{n\in\mathbb Z}\widehat{\drho}_n e^{i\kappa_n x} \quad \text{with} \quad  \widehat{\drho}_0=0. 
\]
Since $
\partial_x\dE= \drho $, we have 
$
\widehat{\dE}_n = \frac{i}{\kappa_n}\widehat{\drho}_n\,,
$
and consequently, 
\begin{align*}
\|(I - P_r)\dE\|_{L_x^2}^2
= \sum_{|k|>r} \frac{1}{|\kappa_k|^2}
   \left|\widehat{\drho}_l\right|^2 
\le \frac{1}{|\kappa_{r+1}|^2}
   \sum_{|k|>r}\left|\widehat{\drho}_k\right|^2 
= \frac{1}{|\kappa_{r+1}|^2}
   \|(I - P_r) \drho\|_{L_x^2}^2.
\end{align*}
If, in addition, $\drho\in H_x^s$, then
\[
\|(I-P_r)\drho\|_{L_x^2}
\le
\inner{\kappa_{r+1}}^{-s}\|\drho\|_{H_x^s},
\]
and consequently
\begin{equation}\label{eq:dE err}
\|(I - P_r)\dE\|_{L_x^2}
\le
|\kappa_{r+1}|^{-1}\inner{\kappa_{r+1}}^{-s}
\norm{\drho}_{H_x^s} \,.
\end{equation}

For the damping part, define
\[
q(x)=\int \delta f(x,v)\,\partial_v\bar f(v)\,dv.
\]
then by Cauchy--Schwarz in $v$, one has 
\[
\|q\|_{H_x^s}
\le
\|\partial_v\bar f\|_{L_v^2}\,\|\delta f\|_{H_x^sL_v^2} \,.
\]
Invoking \eqref{eq:rank truncation err}, we have 
\begin{equation}\label{eq:damping err}
\gamma\|(I-P_r)q\|_{L_x^2}
\le
\gamma\langle \kappa_{r+1}\rangle^{-s}
\|\partial_v\bar f\|_{L_v^2}\,\|\delta f\|_{H_x^sL_v^2}.
\end{equation}
Combining \eqref{eq:dE err} and \eqref{eq:damping err} yields the desired estimate \eqref{eq:finite rank err}.
\end{proof}


\section*{Acknowledgment}
JL and JC were supported by NSF-CCF:2212318, and JC was additionally supported by an Albert and Dorothy Marden Professorship. 
LW is supported in part by NSF grant DMS-1846854, DMS-2513336, and the Simons Fellowship. 

\bibliographystyle{plain}
\bibliography{Ref.bib}
\end{document}